\documentclass[review]{elsarticle}

\usepackage{amsmath}
\usepackage{amssymb}
\usepackage{amsthm}
\usepackage[dvipsnames]{xcolor}
\usepackage{pgfplots}
\usepgfplotslibrary{groupplots}
\usepackage{tikz}
\usetikzlibrary{backgrounds,automata}
\pgfplotsset{width=0.6\textwidth,
    every axis plot/.append style={line width=2pt},
    every axis/.append style={/pgf/number format/.cd,1000 sep={\,},legend pos=outer north east,},
compat=1.3}
\usepackage{caption}
\usepackage{subcaption}
\usepackage{amsmath} 
\usepackage[binary-units]{siunitx}
\usepackage{lineno,hyperref}
\modulolinenumbers[5]
\usepackage[nameinlink]{cleveref}
\usepackage{algorithm}
\usepackage{algpseudocode}

\newtheorem{propostion}{Proposition}

\journal{Journal of \LaTeX\ Templates}

\newcommand{\speedOfLight}{\ensuremath{c}}
\newcommand{\opacityAbsorption}{\ensuremath{\sigma^{\mathrm{a}}}}
\newcommand{\opacityScattering}{\ensuremath{\sigma^{\mathrm{s}}}}
\newcommand{\radiativeConstant}{\ensuremath{a_{\mathrm{r}}}}

\renewcommand{\vector}[1]{\ensuremath{\mathbf{{#1}}}} 
\newcommand{\norm}[1]{\ensuremath{\left \lvert \left \lvert #1 \right \lvert \right \lvert}}
\newcommand{\abs}[1]{\ensuremath{\left \lvert #1 \right \lvert}}
\newcommand{\transpose}[1]{\ensuremath{#1^T}}
\newcommand{\dif}{\ensuremath{\mathrm{d}}}

\newcommand{\radiativeEnergy}{\ensuremath{E_{\mathrm{r}}}}
\newcommand{\radiativeFlux}{\ensuremath{\vector{F}_{\mathrm{r}}}}
\newcommand{\radiativePressure}{\ensuremath{\mathbb{P}_{\mathrm{r}}}}
\newcommand{\reducedFlux}{\ensuremath{f}}
\newcommand{\radiativeTemperature}{\ensuremath{T_{\mathrm{r}}}}
\newcommand{\gasTemperature}{\ensuremath{T_{\mathrm{g}}}}
\newcommand{\gasTemperatureScheme}{\ensuremath{T}}

\newcommand{\hyperbolicVariable}{\ensuremath{\mathcal{U}}}
\newcommand{\hyperbolicFlux}{\ensuremath{\mathcal{F}}}

\newcommand{\radiativeEnergyScheme}{\ensuremath{E}}
\newcommand{\radiativeFluxScheme}{\ensuremath{F}}
\newcommand{\radiativePressureScheme}{\ensuremath{P}}

\newcommand{\spaceStep}{\ensuremath{h}}
\newcommand{\admissibleStatesSet}{\ensuremath{\mathcal{S}}}

\newcommand{\jacobiUnknown}{\ensuremath{v}}
\newcommand{\jacobiRhs}{\ensuremath{b}}
\newcommand{\jacobiResidual}{\ensuremath{r}}
\newcommand{\jacobiError}{\ensuremath{e}}
\newcommand{\jacobiIntermediate}{\ensuremath{u}}
\newcommand{\jacobiOperator}{\ensuremath{\mathcal{A}}}
\newcommand{\leftOperator}{\ensuremath{\mathcal{L}}}
\newcommand{\diagonalOperator}{\ensuremath{\mathcal{D}}}
\newcommand{\rightOperator}{\ensuremath{\mathcal{R}}}
\newcommand{\jacobiCounter}{\ensuremath{k}}

\newcommand{\cfl}{\mathrm{CFL}}
\newcommand{\tolJacobi}{\ensuremath{\varepsilon_J}}
\newcommand{\tolGmg}{\ensuremath{\varepsilon}}
\newcommand{\tolPseudoTimeIncrease}{\ensuremath{\varepsilon_i}}
\newcommand{\tolPseudoTimeDecrease}{\ensuremath{\varepsilon_d}}

\newcommand{\poissonDomain}{\ensuremath{\Omega}}
\newcommand{\poissonMeshFine}{\ensuremath{\poissonDomain_{\spaceStep}}}
\newcommand{\poissonMeshCoarse}{\ensuremath{\poissonDomain_{2 \spaceStep}}}

\newcommand{\restrictionOperator}[3]{\ensuremath{\restrictionOperatorNoArg{#1}{#2} \left( #3 \right)}}
\newcommand{\restrictionOperatorNoArg}[2]{\ensuremath{\mathcal{R}_{#1}^{#2}}}
\newcommand{\prolongationOperator}[3]{\ensuremath{\prolongationOperatorNoArg{#1}{#2} \left( #3 \right)}}
\newcommand{\prolongationOperatorNoArg}[2]{\ensuremath{\mathcal{P}_{#2}^{#1}}}
\newcommand{\fineOperator}{\ensuremath{\jacobiOperator^{\spaceStep}}}
\newcommand{\fineVector}[1]{\ensuremath{\vector{#1}^{\spaceStep}}}
\newcommand{\coarseOperator}{\ensuremath{\jacobiOperator^{2\spaceStep}}}
\newcommand{\coarseVector}[1]{\ensuremath{\vector{#1}^{2\spaceStep}}}
\newcommand{\preSmootherNbIter}{\ensuremath{\smootherNbIter{0}}}
\newcommand{\postSmootherNbIter}{\ensuremath{\smootherNbIter{0}}}
\newcommand{\smootherNbIter}[1]{\ensuremath{\nu_{#1}}}
\newcommand{\indexCellCoarse}[1]{\ensuremath{#1^{2\spaceStep}}}
\newcommand{\indexCellFine}[1]{\ensuremath{#1^{\spaceStep}}}
\newcommand{\indexCellLevel}[2]{\ensuremath{#2^{#1}}}
\newcommand{\nbLevelMax}{\ensuremath{L_{\max}}}
\newcommand{\indexNbCycle}{\ensuremath{\kappa}}
\newcommand{\indexLevel}{\ensuremath{l}}

\newcommand{\cycleIter}[2]{\ensuremath{#1_{(#2)}}} 
\newcommand{\operatorAtLevel}[1]{\ensuremath{\jacobiOperator^{2^{#1}\spaceStep}}} 
\newcommand{\vectorAtLevel}[2]{\vector{#1}^{2^{#2}\spaceStep}} 

\newcommand{\pseudoTime}{\ensuremath{\tau}}
\newcommand{\dTauIm}{\ensuremath{\Delta \pseudoTime^{\text{im}} }}
\newcommand{\dTauEx}{\ensuremath{\Delta \pseudoTime^{\text{ex}} }}
\newcommand{\tauEx}{\ensuremath{\pseudoTime^{\text{ex}} }}
\newcommand{\pseudoTimeProlongation}{\ensuremath{\overline{\tau}}}
\newcommand{\pseudoTimeIndex}{\ensuremath{m}}
\newcommand{\pseudoTimeCounter}{\ensuremath{K}}

\newcommand{\newtonOperator}{\ensuremath{\mathcal{F}}}
\newcommand{\newtonIter}{k}
\newcommand{\newtonUnknown}{\vector{\jacobiUnknown}}
\newcommand{\newtonError}{\delta \newtonUnknown}
\newcommand{\newtonJac}[1]{\frac{\partial \newtonOperator}{\partial \newtonUnknown^{(#1)}} \left( \newtonUnknown^{(#1)} \right)}

\makeatletter
\newenvironment{breakablealgorithm}
{
    \begin{center}
        \refstepcounter{algorithm}
        \hrule height.8pt depth0pt \kern2pt
        \renewcommand{\caption}[2][\relax]{
            {\raggedright\textbf{\fname@algorithm~\thealgorithm} ##2\par}%
            \ifx\relax##1\relax 
                \addcontentsline{loa}{algorithm}{\protect\numberline{\thealgorithm}##2}%
            \else 
                \addcontentsline{loa}{algorithm}{\protect\numberline{\thealgorithm}##1}%
            \fi
            \kern2pt\hrule\kern2pt
        }
    }{
        \kern2pt\hrule\relax
    \end{center}
}
\makeatother

\begin{document}

\begin{frontmatter}

    \title{Towards a multigrid method for the M\textsubscript{1} model for radiative transfer}

    \author[affiliationMdls]{H\'el\`ene Bloch\corref{mycorrespondingauthor}}
    \cortext[mycorrespondingauthor]{Corresponding author}
    \ead{helene.bloch@polytechnique.edu}

    \author[affiliationMdls]{Pascal Tremblin}
    \author[affiliationMatthias]{Matthias Gonz\'alez}
    \author[affiliationMdls]{Edouard Audit}

    \address[affiliationMdls]{Universit\'e Paris-Saclay, UVSQ, CNRS, CEA, Maison de la Simulation, 91191, Gif-sur-Yvette, France}
    \address[affiliationMatthias]{Universit\'e Paris Cit\'e, Universit\'e Paris-Saclay, CEA, CNRS, AIM, F-91191 Gif-sur-Yvette, France}

    \begin{abstract}
        We present a geometric multigrid solver for the M\textsubscript{1} model of radiative transfer without source terms.
        In radiative hydrodynamics applications, the radiative transfer needs to be solved implicitly because of the fast propagation speed of photons relative to the fluid velocity.
        The M\textsubscript{1} model is hyperbolic and can be discretized with an HLL solver, whose time implicit integration can be done using a nonlinear Jacobi method.
        One can show that this iterative method always preserves the admissible states, such as positive radiative energy and reduced flux less than \num{1}.
        To decrease the number of iterations required for the solver to converge, and therefore to decrease the computational cost, we propose a geometric multigrid algorithm.
        Unfortunately, this method is not able to preserve the admissible states.
        In order to preserve the admissible state states, we introduce a pseudo-time such that the solution of the problem on the coarse grid is the steady state of a differential equation in pseudo-time.
        We present preliminary results showing the decrease of the number of iterations and computational cost as a function of the number of multigrid levels used in the method.
        These results suggest that nonlinear multigrid methods can be used as a robust implicit solver for hyperbolic systems such as the M\textsubscript{1} model.
    \end{abstract}


\end{frontmatter}


\section{Introduction}

Radiative transfer describes the transport of radiation, it appears in a wide range of physical situations, from astrophysics, up to medicine.
In radiation therapy, ionizing radiation can be used as a cancer treatment to control malignant cells (e.g., \cite{nascimento2009}).
In astrophysical problems, the radiation exchanges energy with the surrounding gas.
It is the domain of radiation hydrodynamics.
Two phenomena have to be modeled: on one hand, the evolution of the fluid can be described by the Euler equations (e.g., \cite{christodoulou2007}), on the other hand, the radiative transfer is described by an integro-differential equation \cite{chandrasekhar1960,mihalas1984}.
Some source terms model the energy exchange between the fluid and the radiation (e.g., \cite{lowrie1999}).
In this work, we focus on the radiative transfer phenomenon, and we only consider completely transparent media, where the gas and radiation are not coupled.
Photons are in a free-streaming regime, the medium is optically thin.

The quantity of interest for the radiative transfer is the specific intensity $I$.
It depends on six variables: the time $t$, the position $\vector{x}$, the direction of propagation of the photons $\vector{\Omega}$, and the frequency of the photons $\nu$.
Because of this high number of degrees of freedom, solving the integro-differential equation describing the radiative transfer while coupled with a grid-based code for hydrodynamics can be very costly.
Several methods have been developed to tackle this issue, such as ray tracing (e.g., \cite{wise2011}) or Monte-Carlo methods (e.g., \cite{roth2015}).
Ray tracing solves the problem by following the propagation of beams through the fluid.
Although this method is very precise, it is highly costly when coupled with hydrodynamics codes.
The number of degrees of freedom scales with the number of spatial cells multiplied by the number of radiation sources.
Using a Monte-Carlo method, we follow the propagation of ``photon energy packets'' and their interactions with the fluid.
It is based on a stochastic process, making this method accurate, but difficult to couple with a grid based hydrodynamics code.
When the number of photon packets becomes too small in a region, a numerical noise can arise and pollute the simulation results.

In this work, we have chosen to use a gray moment model for the compromise between accuracy and low computational cost.
The specific intensity is averaged over the direction of propagation and the frequency of the photons.
One can consider the first two moments, the radiative energy $\radiativeEnergy$ and the radiative flux $\radiativeFlux$.
A closure relation is needed to link the radiative flux to the radiative energy.
A widely used model is the Flux-Limited Diffusion (FLD) approximation (e.g., \cite{levermore1981}).
However, this model can fail to capture the free-streaming regime of the photons in some cases, see e.g., \cite{gonzalez2007}.
One can consider an extra moment, the radiative pressure $\radiativePressure$.
The closure relation expresses the radiative pressure as a function of radiative energy and radiative flux.
We use the one given by the M\textsubscript{1} model (\cite{levermore1984,dubroca1999}) for its good properties in both optically thin and optically thick media, where the radiation and the fluid strongly interact with each other, even though we do not consider this case.

An explicit solver for the radiative transfer would be restricted by a Courant-Friedrichs-Lewy (CFL) condition, limited by the speed of light.
This will result in a very low time step compared to the hydrodynamics one, which is limited by the speed of sound of the fluid.
We have chosen to use an implicit scheme to solve the M\textsubscript{1} model.
The relation between the radiative energy, radiative flux, and radiative pressure is nonlinear, an implicit solver leads to a nonlinear system.
It can be solved by using a Newton-Raphson method (e.g., \cite{Bertsekas2016}).
At each iteration, a large sparse linear system has to be inverted.
To obtain the best performances,  preconditioned linear solvers have been developed, see for example the discussion in \cite{bloch2021}.
However, large time steps cannot be used in the free-streaming regime because,
in general, the iterations of the Newton-Raphson method do not preserve the admissible states such as positive radiative energy and reduced flux less than \num{1}.

The M\textsubscript{1} model is a hyperbolic system with source terms.
From a numerical point of view, the difficulty comes from a diffusion regime dictated by an asymptotic behavior.
Numerous asymptotic preserving schemes that capture this asymptotic behavior have been developed.
Some of them also preserve the admissible states but are explicit (e.g., \cite{berthon2011,buetDespres2006}), other schemes are implicit but encounter a reduced flux greater than \num{1} especially in the free-streaming regime (e.g., \cite{bloch2021}).
In this work, we only consider the hyperbolic system, without source terms, and we present an implicit scheme that preserves the admissible states.

\cite{pichard2016,pichard2019} suggests solving the nonlinear system with a Jacobi method.
Because this method is iterative, its convergence rate can be improved thanks to multigrid acceleration (e.g., \cite{briggs2000,brandt2011}).
However, this algorithm is not specific to the radiative transfer problem and does not preserve the admissible states.
Inspired by \cite{kifonidis2012}, we introduce a pseudo-time to tackle this issue.

This paper is organized as follows.
We first present in more detail the M\textsubscript{1} model and its discretization in \cref{sect:m1}.
We then show in \cref{sect:gmg:jacobi} the Jacobi method to solve the nonlinear system arising from the discretization of the M\textsubscript{1} model.
We explore the geometric multigrid technique in \cref{sect:gmg:multigrille} and we highlight the modification we made to apply it to radiative transfer.
In \cref{sect:gmg:resultatsNum}, we perform some tests to validate both algorithms.
We also show some performance results.
Finally, we reach our conclusion in \cref{sect:gmg:conc}.

\section{M\textsubscript{1} model and implicit discretization} \label{sect:m1}

\subsection{Continuous model}

The quantity of interest, the specific intensity $I$, describes the rate of radiative transfer of energy, at a point $\vector{x}$ and time $t$.
It verifies the following equation:
\begin{equation}
    \begin{aligned}
        \left( \frac{1}{\speedOfLight} \partial_t + \vector{\Omega} \cdot \nabla \right) I \left( \vector{x}, t, \vector{\Omega}, \nu \right) &= - \left( \opacityAbsorption_{\nu} + \opacityScattering_{\nu} \right) I \left( \vector{x}, t, \vector{\Omega}, \nu \right) + \opacityAbsorption_{\nu} B \left( \nu, \gasTemperature \right) \\
        & + \opacityScattering_{\nu}  \int_{S^2} p^{\nu} \left( \vector{\Omega} \cdot \vector{\Omega}' \right) I \left( \vector{x}, t, \vector{\Omega}', \nu \right) \dif \Omega',
    \end{aligned}
    \label{eq:specificIntensity}
\end{equation}
where $\speedOfLight$ is the speed of light, $\nu$ is the frequency of the photons.
$\vector{\Omega}$ is the direction of propagation of the photons: in the vacuum, photons propagate in a straight line at velocity $\speedOfLight \vector{\Omega}$.
$\opacityAbsorption_{\nu}$ is the absorption coefficient, $\opacityScattering_{\nu}$ is the scattering coefficient, and $B$ is the black body specific intensity.
See, for example, \cite{chandrasekhar1960} for the derivation of \cref{eq:specificIntensity}.

Let us now define the gray three first moments of the specific intensity as
\begin{equation}
    \begin{aligned}
        \radiativeEnergy &= \frac{1}{\speedOfLight} \int_0^{\infty} \int_{S^2} I \left( \vector{x}, t, \vector{\Omega}, \nu \right) \dif \Omega \dif \nu\\
        \radiativeFlux &= \int_0^{\infty} \int_{S^2} I \left( \vector{x}, t, \vector{\Omega}, \nu \right) \vector{\Omega} \dif \Omega \dif \nu\\
        \radiativePressure &= \frac{1}{\speedOfLight} \int_0^{\infty} \int_{S^2} I \left( \vector{x}, t, \vector{\Omega}, \nu \right) \vector{\Omega} \otimes \vector{\Omega} \dif \Omega \dif \nu.\\
    \end{aligned}
    \label{eq:defMoments}
\end{equation}
$\radiativeEnergy$ is the radiative energy, $\radiativeFlux$ is the radiative flux, and $\radiativePressure$ is the radiative pressure.

Using \cref{eq:specificIntensity}, one can show that the moments obey the following system:
\begin{equation}
    \begin{aligned}
        \partial_t \radiativeEnergy + \nabla \cdot \radiativeFlux &= \speedOfLight \left( \sigma_{\mathrm{P}} \radiativeConstant \gasTemperature^4 - \sigma_{\mathrm{E}} \radiativeEnergy \right)\\
        \partial_t \radiativeFlux + \speedOfLight^2 \nabla \cdot \radiativePressure &= - \speedOfLight \left( \sigma_{\mathrm{F}} + \left( 1 - g \right) \sigma_{\mathrm{R}} \right) \radiativeFlux,
    \end{aligned}
    \label{eq:m1}
\end{equation}
with
$\sigma_{\mathrm{E}}$ (resp. $\sigma_{\mathrm{P}}$) the mean over frequency of $\opacityAbsorption_{\nu}$ weighted by the radiative energy (resp. the black body specific intensity), $\sigma_{\mathrm{F}}$ (resp. $\sigma_{\mathrm{S}}$) the mean of $\opacityAbsorption_{\nu}$ (resp. $\opacityScattering_{\nu}$) weighted by the radiative flux, $g$ the first moment of the phase function $p^{\nu}$,
$\gasTemperature$ the temperature of the gas, and $\radiativeConstant$ the radiative constant.
See \cite{gonzalez2007} for the discussion about the source terms and the approximations that can be made, such as $\sigma_{\mathrm{E}} = \sigma_{\mathrm{P}}$.

Using \cref{eq:defMoments}, one can notice that $\radiativeEnergy > 0$ because $I > 0$, and 
\begin{equation}
    \norm{\radiativeFlux} \le \int_0^{\infty} \int_{S^2} \underbrace{I \left( \vector{x}, t, \vector{\Omega}, \nu \right)}_{\ge 0} \underbrace{\norm{\vector{\Omega}}}_{= 1} \dif \Omega \dif \nu = \speedOfLight \radiativeEnergy.
\end{equation}
This condition rewrites $\reducedFlux \le 1$, where $\reducedFlux = \norm{\vector{\reducedFlux}} = \frac{\norm{\vector{\radiativeFlux}}}{\speedOfLight \radiativeEnergy}$ is the reduced flux. It ensures that the radiative energy cannot be transported faster than the speed of light.
These two conditions, the positive radiative energy and the reduced flux less than \num{1}, have to be preserved by the numerical scheme, for physical and numerical reasons.

We now specify the closure relation, which is a way to express $\radiativePressure$ as a function of $\radiativeEnergy$ and $\radiativeFlux$.
We have chosen to use the M\textsubscript{1} closure relation for its good properties in the free-streaming regime.
We have
\begin{equation}
    \radiativePressure = \left( \frac{1-\chi}{2} \mathbb{I} + \frac{3 \chi - 1}{2} \vector{n} \otimes \vector{n} \right) \radiativeEnergy,
    \label{eq:closureRelation}
\end{equation}
where $\vector{n} = \frac{\vector{\reducedFlux}}{\reducedFlux}$ is a unit vector aligned with the reduced flux, $\reducedFlux = \norm{\vector{\reducedFlux}}$, and $\mathbb{I}$ is the identity matrix.
$\chi$ is the Eddington factor, defined as $\chi = \frac{3+4f^2}{5+2\sqrt{4-3\reducedFlux^2}}$.
\Cref{eq:closureRelation} can be obtained by applying a Lorentz transform to an isotropic distribution of photons (\cite{levermore1984}) or by maximizing the radiative entropy (\cite{dubroca1999}).

More details about the derivation of the model and the closure relation can be found in \cite{these}.

Developing numerical schemes to solve \cref{eq:m1} can be challenging due to the presence of the source terms and their behavior in optically thick media.
Several numerical schemes that take into account the source terms have been developed, see e.g., \cite{jin1996,buetDespres2008,berthon2011}.
However, most of them are explicit schemes.
For some physical applications, an implicit scheme is mandatory.
For example, in radiation hydrodynamics, when radiative transfer is coupled with hydrodynamics, the speed of light is larger than the speed of sound of the fluid.
A solution to avoid a time step restricted by the speed of light for the hydrodynamics is to use an implicit solver for the radiative transfer.
When both phenomena are coupled, source terms are mandatory to model this coupling.
Developing an implicit scheme with source term is still an active topic of research.
Here, we do not address the source terms, but we focus on the implicit property of the scheme
In the following, we set
$\sigma_{\mathrm{E}} = \sigma_{\mathrm{P}} = \sigma_{\mathrm{F}} = \sigma_{\mathrm{S}} = 0$
which allows us to use standard numerical schemes for \cref{eq:m1}.

One can show that the M\textsubscript{1} model is hyperbolic and can be written as
\begin{equation*}
  \partial_t \hyperbolicVariable + \nabla \cdot \hyperbolicFlux \left( \hyperbolicVariable \right) = 0,
\end{equation*}
with $\hyperbolicVariable = \left( \radiativeEnergy, \vector{\radiativeFlux}  \right)$ and $\hyperbolicFlux$ the flux function, defined as $\hyperbolicFlux \left( \hyperbolicVariable \right) = \left( \vector{\radiativeFlux}, \speedOfLight^2 \radiativePressure \right)$.
The discretization of \cref{eq:m1} is discussed in the next section, without considering the source terms anymore.

\subsection{HLL solver}

Let us first introduce some notations.
We note $\spaceStep$ the step along the x-direction.
$\Delta t$ is the time interval between the current time $t^n$ and $t^{n+1}$.
We write $x_{i}$ the center of the cell $i$.
We use the notation $u_{i}^n$ to represent the averaged quantity associated with the field $u$ at time $t^n$ in cell $i$ (finite volume). 
To ease notation, we drop the indices $\mathrm{r}$ for all radiative variables

The M\textsubscript{1} model is hyperbolic, we can discretize it using a time implicit HLL solver
\cite{harten1983}.
For the sake of simplicity, we only consider the one-dimensional case, but the extension to higher dimensions is straightforward.
The flux at the interface $i + \frac{1}{2}$ between cells $i$ and $i+1$ depends on the speed of propagation, $\lambda_{i+\frac{1}{2}}^+$ and $\lambda_{i+\frac{1}{2}}^-$.
$\lambda_{i+\frac{1}{2}}^+$ and $\lambda_{i+\frac{1}{2}}^-$ can be computed using the Eddington factor $\chi$.
One can show that $\lambda_{i+\frac{1}{2}}^+ \le \speedOfLight$ and $\lambda_{i+\frac{1}{2}}^- \ge - \speedOfLight$.
Using these upper bounds is sufficient for the scheme.
Therefore, still for the sake of simplicity, we use these upper bounds in the flux.
One can see the discussion in \cite{gonzalez2007} about the impact of this choice.
This is equivalent to use the Rusanov's flux \cite{rusanov1962}.
This leads to solving the following nonlinear system:
\begin{equation}
    \begin{aligned}
        \radiativeEnergyScheme^{n+1}_{i} &\left( 1 + \speedOfLight \frac{\Delta t}{\spaceStep}  \right) - \frac{\Delta t}{2 \spaceStep} \left( \speedOfLight \radiativeEnergyScheme^{n+1}_{i+1} - \radiativeFluxScheme^{n+1}_{i+1} \right) - \frac{\Delta t}{2 \spaceStep} \left( \speedOfLight \radiativeEnergyScheme^{n+1}_{i-1} + \radiativeFluxScheme^{n+1}_{i-1} \right) &= \radiativeEnergyScheme^n_{i} \\
        \radiativeFluxScheme^{n+1}_{i} &\left( 1 + \speedOfLight \frac{\Delta t}{\spaceStep} \right) - \frac{\speedOfLight \Delta t}{2 \spaceStep} \left( \radiativeFluxScheme^{n+1}_{i+1} - \speedOfLight \radiativePressureScheme^{n+1}_{i+1} \right) - \frac{\speedOfLight \Delta t}{2 \spaceStep} \left( \radiativeFluxScheme^{n+1}_{i-1} + \speedOfLight \radiativePressureScheme^{n+1}_{i-1} \right) &= \radiativeFluxScheme^n_{i}.\\
    \end{aligned}
    \label{eq:gmg:jacobi:hll}
\end{equation}
\Cref{eq:gmg:jacobi:hll} can be solved using a Newton-Raphson method, as done in \cite{bloch2021,gonzalez2007} for example.
However, numerical tests close to the free-streaming limit have shown that the admissible states are not preserved when large time steps are used.
A solution is therefore to reduce the time step, but it leads to poor performances when the radiative transfer is coupled to hydrodynamics.

In the next section, we present another method to solve \cref{eq:gmg:jacobi:hll} while preserving the admissible states.

\section{Nonlinear Jacobi method} \label{sect:gmg:jacobi}

\subsection{Algorithm}

Let us follow the work of \cite{pichard2016,pichard2019}.
We first define the set of admissible states
\begin{equation}
    \admissibleStatesSet = \left\{ \left( \radiativeEnergy, \radiativeFlux \right) , \radiativeEnergy > 0, \norm{\radiativeFlux} \le \speedOfLight \radiativeEnergy \right\}.
    \label{eq:admissibleStateSet}
\end{equation}
By writing $\vector{\jacobiUnknown} = \transpose{\left( \cdots, \radiativeEnergyScheme_i, \vector{\radiativeFluxScheme}_i, \cdots \right)} \in \admissibleStatesSet^{N}$, where $N$ is the number of cells, \cref{eq:gmg:jacobi:hll} rewrites 
\begin{equation}
  - \leftOperator^h \left( \jacobiUnknown^{n+1}_{i-1} \right) + \diagonalOperator^h \left( \jacobiUnknown^{n+1}_{i} \right) - \rightOperator^h \left( \jacobiUnknown^{n+1}_{i+1} \right) = \jacobiUnknown^{n}_{i},
    \label{eq:gmg:jacobi:systSolved}
\end{equation}
where $\jacobiUnknown_i = \left( \radiativeEnergyScheme_i, \vector{\radiativeFluxScheme}_i \right)$.
The operators $\leftOperator^h$, $\diagonalOperator^h$, and $\rightOperator^h$ contain the terms depending on $\left( \radiativeEnergyScheme^{n+1}_{i-1}, \vector{\radiativeFluxScheme}^{n+1}_{i-1} \right)$, $\left( \radiativeEnergyScheme^{n+1}_{i}, \vector{\radiativeFluxScheme}^{n+1}_{i} \right)$ and $\left( \radiativeEnergyScheme^{n+1}_{i+1}, \vector{\radiativeFluxScheme}^{n+1}_{i+1} \right)$ respectively.
In 1D,
\begin{equation}
  \leftOperator^h\left( \jacobiUnknown \right) = 
  \begin{pmatrix}
    \frac{\Delta t}{2 \spaceStep} \left( \speedOfLight \radiativeEnergyScheme + \radiativeFluxScheme \right)\\
    \frac{\speedOfLight \Delta t}{2 \spaceStep} \left( \radiativeFluxScheme + \speedOfLight \radiativePressureScheme \right)
  \end{pmatrix}
  , \ 
  \diagonalOperator^h \left( \jacobiUnknown \right) =
  \begin{pmatrix}
    \left( 1 + \frac{\speedOfLight \Delta t}{\spaceStep} \right) \radiativeEnergyScheme \\
    \left( 1 + \frac{\speedOfLight \Delta t}{\spaceStep} \right) \radiativeFluxScheme
  \end{pmatrix}
  , \ 
  \rightOperator^h\left( \jacobiUnknown \right) = 
  \begin{pmatrix}
    \frac{\Delta t}{2 \spaceStep} \left( \speedOfLight \radiativeEnergyScheme - \radiativeFluxScheme \right)\\
    \frac{\speedOfLight \Delta t}{2 \spaceStep} \left( \radiativeFluxScheme - \speedOfLight \radiativePressureScheme \right)
  \end{pmatrix}
  .
\end{equation}
To ease notations, we write $\leftOperator = \leftOperator^h$, $\diagonalOperator = \diagonalOperator^h$, and $\rightOperator = \rightOperator^h$ when there is no ambiguity.

From \cite{pichard2016}, we solve \cref{eq:gmg:jacobi:systSolved} using \cref{algo:jacobi}.

\begin{algorithm}
    \caption{Nonlinear Jacobi method}
    \begin{algorithmic}[1]
        \State Initialization: $\vector{\jacobiUnknown}^{n+1, (0)} = \vector{\jacobiUnknown}^n$
        \While{$\norm{\jacobiOperator \left( \vector{\jacobiUnknown}^{n+1, (\jacobiCounter)} \right) - \vector{\jacobiRhs}} > \tolJacobi$}
        \For{each cell $i$}
        \State 
        \begin{equation}
            \jacobiUnknown^{n+1, (\jacobiCounter+1)}_{i} = \diagonalOperator^{-1} \left( \jacobiRhs_{i} + \leftOperator \left( \jacobiUnknown^{n+1, (\jacobiCounter)}_{i-1} \right) + \rightOperator \left( \jacobiUnknown^{n+1, (\jacobiCounter)}_{i+1} \right) \right)
            \label{eq:gmg:jacobi:algo}
        \end{equation}
        \EndFor
        \State $\jacobiCounter \gets \jacobiCounter +1$
        \EndWhile
        \State $\vector{\jacobiUnknown}^{n+1} = \vector{\jacobiUnknown}^{n+1, (\jacobiCounter)}$
    \end{algorithmic}
    \label{algo:jacobi}
\end{algorithm}

\Cref{eq:gmg:jacobi:hll} can also be seen as 
\begin{equation}
  \jacobiOperator^h \left( \vector{\jacobiUnknown} \right) = \vector{\jacobiRhs},
  \label{eq:gmg:nonLinearSyst}
\end{equation}
where $\jacobiOperator^h$ is a nonlinear operator, $\vector{\jacobiUnknown}$ is the vector of unknowns, and $\vector{\jacobiRhs}$ is a known vector.
As before, we drop the $h$ exponent whenever possible.
The $i$-th row of $\jacobiOperator$ is
\begin{equation}
  \left( \jacobiOperator \left( \vector{\jacobiUnknown} \right) \right)_i = 
  \begin{pmatrix}
    \radiativeEnergyScheme^{n+1}_{i} \left( 1 + \speedOfLight \frac{\Delta t}{\spaceStep}  \right) - \frac{\Delta t}{2 \spaceStep} \left( \speedOfLight \radiativeEnergyScheme^{n+1}_{i+1} - \radiativeFluxScheme^{n+1}_{i+1} \right) - \frac{\Delta t}{2 \spaceStep} \left( \speedOfLight \radiativeEnergyScheme^{n+1}_{i-1} + \radiativeFluxScheme^{n+1}_{i-1} \right) \\
    \radiativeFluxScheme^{n+1}_{i} \left( 1 + \speedOfLight \frac{\Delta t}{\spaceStep} \right) - \frac{\speedOfLight \Delta t}{2 \spaceStep} \left( \radiativeFluxScheme^{n+1}_{i+1} - \speedOfLight \radiativePressureScheme^{n+1}_{i+1} \right) - \frac{\speedOfLight \Delta t}{2 \spaceStep} \left( \radiativeFluxScheme^{n+1}_{i-1} + \speedOfLight \radiativePressureScheme^{n+1}_{i-1} \right)
  \end{pmatrix}
  .
\end{equation}
Let us notice that if $\jacobiOperator$ is a linear operator, \cref{algo:jacobi} simplifies into the classical Jacobi  method for a tridiagonal matrix (see e.g., \cite{saad2003}).

It is shown in \cite{pichard2016} that this algorithm converges to the unique solution $\vector{\jacobiUnknown}^{n+1}$ because the operator $\jacobiOperator$ is contractant.
The proof relies on the study of the eigenvalues of the Jacobian of $\jacobiOperator$.

\subsection{Preservation of the admissible states} \label{sect:jacobi:etatsAdmissibles}

\begin{propostion}
  If $\jacobiUnknown^n_i \in \admissibleStatesSet$, then $\jacobiUnknown^{n+1}_i$ obtained with \cref{algo:jacobi} is also admissible, i.e., $\jacobiUnknown^{n+1}_i \in \admissibleStatesSet$.
\end{propostion}

\begin{proof}
  Let us show that if $\vector{\jacobiUnknown}^{n+1, (\jacobiCounter)} \in \admissibleStatesSet^N$, then $\vector{\jacobiUnknown}^{n+1, (\jacobiCounter+1)} \in \admissibleStatesSet^N$.

  Let us recall Proposition~5.1 in \cite{pichard2016}. 
  We write it with our notations, but it can be used in a more general case.
  To fit \cite{pichard2016}'s notations, we introduce $\vector{m} = \left(\frac{1}{\speedOfLight}, \vector{\Omega} \right)$.

  ``Suppose $\hyperbolicVariable \in \admissibleStatesSet$ and the flux $\hyperbolicFlux \left( \hyperbolicVariable \right)$ is defined from a realizable closure, i.e. such that 
  \begin{equation*}
    \exists I > 0 \text{ s.t. } \hyperbolicVariable = \int_0^{\infty} \int_{S^2} \vector{m} I \dif \Omega \dif \nu, \ \hyperbolicFlux \left( \hyperbolicVariable \right) = \int_0^{ \infty} \int_{S^2} \speedOfLight \vector{\Omega} \vector{m} I \dif \Omega \dif \nu,
  \end{equation*}
  then
  \begin{equation*}
    \hyperbolicVariable \pm \hyperbolicFlux \left( \hyperbolicVariable \right) \in \admissibleStatesSet."
  \end{equation*}

        Furthermore, let us notice that the sum of two elements in $\admissibleStatesSet$ and the multiplication of an element of $\admissibleStatesSet$ by a positive scalar are also in $\admissibleStatesSet$ (but $\admissibleStatesSet$ is not a vector space because an element of $\admissibleStatesSet$ multiplied by a negative scalar is not in $\admissibleStatesSet$).
        One can then show that $\leftOperator$ and $\rightOperator$ are stable, i.e., $\leftOperator \left( \jacobiUnknown \right) \in \admissibleStatesSet$ and $\rightOperator \left( \jacobiUnknown \right) \in \admissibleStatesSet$ if $\jacobiUnknown \in \admissibleStatesSet$.

    Therefore, 
    because $\jacobiRhs_i = \jacobiUnknown^n_i$, it is admissible and $\jacobiRhs_i + \leftOperator \left( \jacobiUnknown_{i-1} \right) + \rightOperator \left( \jacobiUnknown_{i+1} \right) \in \admissibleStatesSet$.

    Finally, in this particular case, $\diagonalOperator = \left( 1 + \speedOfLight \frac{\Delta t}{\spaceStep} \right) \mathbb{I}$, where $\mathbb{I}$ is the identity matrix.
    So, if $\jacobiUnknown \in \admissibleStatesSet$, then $\diagonalOperator^{-1} \left( \jacobiUnknown \right) = \frac{1}{1 + \speedOfLight \frac{\Delta t}{\spaceStep}} \jacobiUnknown \in \admissibleStatesSet$.

    We have shown that $\diagonalOperator^{-1}$, $\leftOperator$, and $\rightOperator$ are stable.
    Therefore, $\jacobiUnknown^{n+1, (\jacobiCounter+1)}_{i}$ is admissible.
    By induction, the proof is complete.

\end{proof}

\section{Geometric multigrid (GMG)} \label{sect:gmg:multigrille}

\subsection{Convergence problematic} \label{sect:gmg:interet}

\begin{figure}
    \centering
    \begin{tikzpicture}
  \begin{semilogyaxis}[
      xlabel = {Number of iterations ($\jacobiCounter$)},
      ylabel = {Residual $\frac{\norm{\vector{\jacobiResidual}^{(\jacobiCounter)}} }{ \norm{\vector{\jacobiResidual}^{(0)}} }$},
      axis lines = left,
    ]

    \addplot[color=red] table {images/iterJacobiCv/output_jacobi_129.txt};
    \addlegendentry{$129\times129$ cells};

    \addplot[color=green] table {images/iterJacobiCv/output_jacobi_257.txt};
    \addlegendentry{$257\times257$ cells};

    \addplot[color=blue] table {images/iterJacobiCv/output_jacobi_513.txt};
    \addlegendentry{$513\times513$ cells};

  \end{semilogyaxis}
\end{tikzpicture}
    \caption[Residual as a function of the number of iterations for Jacobi method]{Evolution of the residual as a function of the number of iterations of Jacobi method, with different resolutions.}
    \label{fig:gmg:cv:iterJacobi}
\end{figure}
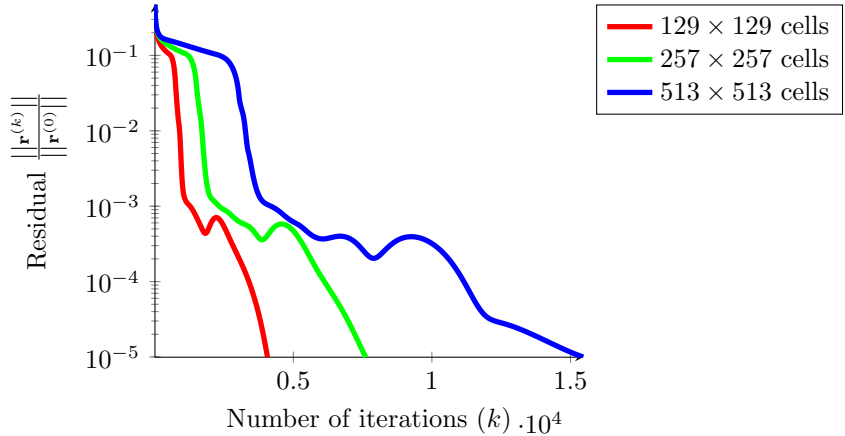

We investigate the performance of the nonlinear Jacobi algorithm on an academic test.
The one considered here is the beam test described in \cite{gonzalez2007,richling2001}.
It is the propagation of a beam in the free-streaming regime.
The domain is $[-1, 1] \times [-1,1]$ and it is discretized with the same number of cells in the $x$ and $y$-directions.
The initial temperature is $\gasTemperatureScheme_0 = \radiativeTemperature = \SI{300}{\kelvin}$, where $\radiativeTemperature = \left( \frac{\radiativeEnergy}{\radiativeConstant} \right)^{\frac{1}{4}}$ is the radiative temperature.
The initial radiative flux is $\radiativeFlux = \vector{0}$.
At time $t = 0$, a beam is introduced with $\radiativeTemperature = \SI{1000}{\kelvin}$
and $\reducedFlux = 1$
with an angle of $45^\circ$ at $x = -1$ and $y \in [-0.875,-0.75]$.
Because the solver is time-implicit, we set a time step $\Delta t$ such that the steady state is reached with one time iteration.

\begin{table}
    \centering
    \begin{tabular}{cc}
        \hline
        Number of cells & Cell-updates/\si{\second} \\
        \hline
        $129 \times 129$ & \num{146}\\
        $257 \times 257$ & \num{76}\\
        $513 \times 513$ & \num{38}
    \end{tabular}
    \caption[Number of cell-updates per second of Jacobi method]{Number of cell-updates per second of Jacobi method, with different resolutions.}
    \label{tab:gmg:cv:perfJacobi}
\end{table}

Let us define the residual as
\begin{equation}
    \vector{\jacobiResidual}^{(\jacobiCounter)} = \vector{\jacobiRhs} - \jacobiOperator \left( \vector{\jacobiUnknown}^{(\jacobiCounter)} \right).
    \label{eq:defResidual}
\end{equation}
\Cref{fig:gmg:cv:iterJacobi} shows the evolution of the norm of the residual as a function of the number of iterations of the Jacobi method, with different resolutions. 
If the operator $\jacobiOperator$ were linear, one could expect the residual to decrease linearly with the number of iterations, see, for example, results in \cite{gander2010} obtained with a linear problem.
As the resolution increases, the number of iterations needed to reach the same residual also increases, from \num{4000} iterations with $129 \times 129$ cells up to \num{15000} iterations with $513 \times 513$ cells. 
As shown by \cref{tab:gmg:cv:perfJacobi}, performances decrease as the resolution increases, from around \num{150} cell-updates/\si{\second} with $129 \times 129$ cells down to \num{40} cell-updates/\si{\second} with $513 \times 513$ cells.
If the number of iterations of the Jacobi method were constant with the resolution, the number of cell-updates/\si{\second} would be constant.
Because the method requires more iterations to reach the same residual when the resolution increases, performances decrease.

Indeed, \cite{saad2003} explains that the convergence rate decreases as the size of the problem increases, resulting in a decrease in performance.
For linear problems, classical iterative methods are efficient to compute the high frequencies of the solution, but lack efficiency to compute its low frequencies.
However, the computation is easier on a coarser grid with fewer unknowns.
This observation has led to the development of the geometric multigrid (GMG) technique. 
The initial guess for the Jacobi algorithm (or any iterative method) is an interpolation of a solution computed on a coarser grid.
This method requires information about the geometry of the problem, unlike preconditioners based on Krylov subspace. 
This additional information allows the multigrid method to be very efficient, but it lacks generality.
To tackle this issue, algebraic multigrid methods have also been developed, see e.g., \cite{saad2003}.

In this work, we focus on the GMG technique, and we first present the Full Approximation Scheme (FAS) in \cref{sect:gmg:fas} that can be used to solve a nonlinear system arising from the discretization of a differential equation.
Then, in \cref{sect:gmg:m1} we apply it to the M\textsubscript{1} model, and we highlight the distinctive features of radiative transfer, i.e., the preservation of the admissible states $\radiativeEnergy > 0$ and $\reducedFlux \le 1$. 

\subsection{General case} \label{sect:gmg:fas}

In this section, we present the main ideas of the geometric multigrid method.
It does not intend to be a full review of existing work, see e.g., \cite{briggs2000,brandt2011}.

We are interested in solving the nonlinear system \cref{eq:gmg:nonLinearSyst} with the Jacobi method presented in \cref{sect:gmg:jacobi}.
We assume that this system is obtained by the discretization of a differential equation on the domain $\poissonDomain$ and therefore depends on a step $\spaceStep$.

As mentioned above, the initial guess used in the Jacobi algorithm is the interpolation of a solution obtained on a coarser grid.
This solution computed on the coarse grid can be itself computed thanks to a solution obtained on a third grid, even coarser.
This process can be applied recursively until the last grid has only a few unknowns and the problem could be solved with a direct method, which leads to the so-called ``nested iterations''.
From \cite{briggs2000}, this requires the definition of two elements:
\begin{itemize}
    \item a hierarchy of grids, completed with restriction and prolongation operators;
    \item an iterative method used as a smoother to solve the system on a given mesh.
\end{itemize}
In the linear case, one can study the eigenvalues and eigenvectors of $\jacobiOperator$, as well as the iterative method to understand how the smoother allows the quick decrease of the high frequencies of the error.

The hierarchy of grids is handled through recursiveness, we only need to describe the method with two grids.
The domain $\poissonDomain$ is discretized with two Cartesian meshes: the first one with step $\spaceStep$, written $\poissonMeshFine$ and the second one with step $2\spaceStep$, written $\poissonMeshCoarse$.
$\poissonMeshFine$ contains $2^d$ times more elements than $\poissonMeshCoarse$, where $d$ is the dimension of the problem.

Each vector or operator is defined on a mesh.
To highlight the dependence on the mesh, we write all the objects with the size of the grid it is defined on.
For example, we write now $\fineVector{\jacobiUnknown}$ instead of $\vector{\jacobiUnknown}$.

Let us now define the restriction and prolongation operators used for inter-grid operations.
We write the prolongation operator $\prolongationOperatorNoArg{\spaceStep}{2\spaceStep} : \poissonMeshCoarse \rightarrow \poissonMeshFine$ and the restriction operator $\restrictionOperatorNoArg{\spaceStep}{2\spaceStep} : \poissonMeshFine \rightarrow \poissonMeshCoarse$.

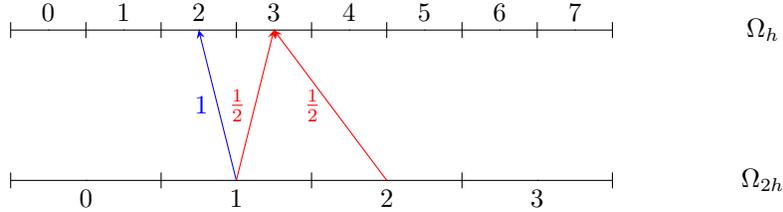
\begin{figure}
    \centering
    \newcommand{\highFineMesh}{2}
\newcommand{\coeffLabels}{\highFineMesh cm}

\begin{tikzpicture}
  \draw (0,\highFineMesh) -- (8,\highFineMesh);
  \foreach \x in {0,...,8} \draw (\x,\coeffLabels+3pt) -- (\x,\coeffLabels-3pt) node {};
  \foreach \x in {0,...,7} \draw (\x+0.5, \coeffLabels) -- (\x+0.5,\coeffLabels) node[above] {\x};
  \node at (10,\highFineMesh) {$\poissonMeshFine$};

  \draw (0,0) -- (8,0);
  \foreach \x in {0,...,4} \draw (2*\x,3pt) -- (2*\x,-3pt) node {};
  \foreach \x in {0,...,3} \draw (2*\x+1,0) -- (2*\x+1,0) node[below] {\x};
  \node at (10,0) {$\poissonMeshCoarse$};

  \draw[color=blue,-stealth]  (3,0) -- (2.5,\highFineMesh);
  \node[color=blue, left] at (2.75, \highFineMesh/2) {$1$};

  \draw[color=red,-stealth] (3,0) -- (3.5,\highFineMesh);
  \node[color=red, left] at (3.25, \highFineMesh/2) {$\frac{1}{2}$};
  \draw[color=red,-stealth] (5,0) -- (3.5,\highFineMesh);
  \node[color=red, left] at (4.25, \highFineMesh/2) {$\frac{1}{2}$};
\end{tikzpicture}
    \caption[Prolongation operator in the one-dimensional case]{Prolongation operator in the one-dimensional case
      with weights.
      Ticks represent the interfaces between two cells, and the indices of the cells are located at the center.
      The blue arrow represents the operator if $\indexCellFine{i}$ is even, and the red arrows represent the operator if $\indexCellFine{i}$ is odd.
    }
    \label{fig:prolongation1D}
\end{figure}

We use full weighting operators (\cite{strang2006}).
Let us begin with the prolongation operator.
Let us consider the one-dimensional case with an 
even
number of cells in the fine mesh.
Cells with index $\indexCellFine{i}$ even in the fine mesh are the same as cells in the coarse mesh with index $\indexCellCoarse{i}$
(see \cref{fig:prolongation1D}). 
Values corresponding to these cells in the coarse mesh are just moved in the fine mesh.
The other values in the fine mesh are obtained by linear interpolation.
The prolongation operator writes
\begin{equation}
    \left( \prolongationOperator{\spaceStep}{2\spaceStep}{\coarseVector{\jacobiUnknown}} \right)_{\indexCellFine{i}} = 
    \left\{
        \begin{aligned}
            &\coarseVector{\jacobiUnknown}_{\indexCellCoarse{i}} \text{ if } \indexCellFine{i} \text{ is even},\\
            &\frac{1}{2} \left( \coarseVector{\jacobiUnknown}_{\indexCellCoarse{i}} + \coarseVector{\jacobiUnknown}_{\indexCellCoarse{i}+1} \right) \text{ if } \indexCellFine{i} \text{ is odd}.
        \end{aligned}
    \right.
    \label{eq:prolongationOperator1D}
\end{equation}

\begin{figure}
  \centering
  \newcommand{\multiplicationFactor}{0.4}

\begin{tikzpicture}
  \draw (\multiplicationFactor*0, \multiplicationFactor*10) -- (\multiplicationFactor*10, \multiplicationFactor*10);
  \draw [dashed] (\multiplicationFactor*1.5, \multiplicationFactor*12) -- (\multiplicationFactor*11.5, \multiplicationFactor*12);
  \draw (\multiplicationFactor*3, \multiplicationFactor*14) -- (\multiplicationFactor*13, \multiplicationFactor*14);
  \draw [dashed] (\multiplicationFactor*4.5, \multiplicationFactor*16) -- (\multiplicationFactor*14.5, \multiplicationFactor*16);
  \draw (\multiplicationFactor*6, \multiplicationFactor*18) -- (\multiplicationFactor*16, \multiplicationFactor*18);

  \draw (\multiplicationFactor*0, \multiplicationFactor*10) -- (\multiplicationFactor*6, \multiplicationFactor*18);
  \draw [dashed] (\multiplicationFactor*2.5, \multiplicationFactor*10) -- (\multiplicationFactor*8.5, \multiplicationFactor*18);
  \draw (\multiplicationFactor*5, \multiplicationFactor*10) -- (\multiplicationFactor*11, \multiplicationFactor*18);
  \draw [dashed] (\multiplicationFactor*7.5, \multiplicationFactor*10) -- (\multiplicationFactor*13.5, \multiplicationFactor*18);
  \draw (\multiplicationFactor*10, \multiplicationFactor*10) -- (\multiplicationFactor*16, \multiplicationFactor*18);

  \node [left] at (\multiplicationFactor*2, \multiplicationFactor*11) {$\indexCellFine{i}, \indexCellFine{j}$};

  \draw (\multiplicationFactor*0, \multiplicationFactor*0) -- (\multiplicationFactor*10, \multiplicationFactor*0);
  \draw (\multiplicationFactor*3, \multiplicationFactor*4) -- (\multiplicationFactor*13, \multiplicationFactor*4);
  \draw (\multiplicationFactor*6, \multiplicationFactor*8) -- (\multiplicationFactor*16, \multiplicationFactor*8);

  \draw (\multiplicationFactor*0, \multiplicationFactor*0) -- (\multiplicationFactor*6, \multiplicationFactor*8);
  \draw (\multiplicationFactor*5, \multiplicationFactor*0) -- (\multiplicationFactor*11, \multiplicationFactor*8);
  \draw (\multiplicationFactor*10, \multiplicationFactor*0) -- (\multiplicationFactor*16, \multiplicationFactor*8);

  \node [below] at (\multiplicationFactor*4, \multiplicationFactor*2) {$\indexCellCoarse{i}, \indexCellCoarse{j}$};
  \node [below] at (\multiplicationFactor*9, \multiplicationFactor*2) {$\indexCellCoarse{i}+1, \indexCellCoarse{j}$};
  \node [below] at (\multiplicationFactor*7, \multiplicationFactor*6) {$\indexCellCoarse{i}, \indexCellCoarse{j}+1$};
  \node [right] at (\multiplicationFactor*12, \multiplicationFactor*6) {$\indexCellCoarse{i}+1, \indexCellCoarse{j}+1$};

  \draw [color=blue, -stealth] (\multiplicationFactor*4, \multiplicationFactor*2) -- (\multiplicationFactor*2, \multiplicationFactor*11);
  \node[color=blue, left] at (\multiplicationFactor*3, \multiplicationFactor*6.5) {$1$};

  \draw [color=red, -stealth] (\multiplicationFactor*4, \multiplicationFactor*2) -- (\multiplicationFactor*4.4, \multiplicationFactor*11);
  \draw [color=red, -stealth] (\multiplicationFactor*9, \multiplicationFactor*2) -- (\multiplicationFactor*4.6, \multiplicationFactor*11);
  \node[color=red] at (\multiplicationFactor*4.5, \multiplicationFactor*7.5) {$\frac{1}{2}$};

  \draw [color=green, -stealth] (\multiplicationFactor*4, \multiplicationFactor*2) -- (\multiplicationFactor*3.4, \multiplicationFactor*13);
  \draw [color=green, -stealth] (\multiplicationFactor*7, \multiplicationFactor*6) -- (\multiplicationFactor*3.6, \multiplicationFactor*13);
  \node[color=green, left] at (\multiplicationFactor*3.695, \multiplicationFactor*7.5) {$\frac{1}{2}$};

  \draw [color=Orange, -stealth] (\multiplicationFactor*4, \multiplicationFactor*2) -- (\multiplicationFactor*5.9, \multiplicationFactor*12.9);
  \draw [color=Orange, -stealth] (\multiplicationFactor*9, \multiplicationFactor*2) -- (\multiplicationFactor*6.1, \multiplicationFactor*12.9);
  \draw [color=Orange, -stealth] (\multiplicationFactor*7, \multiplicationFactor*6) -- (\multiplicationFactor*5.9, \multiplicationFactor*13.1);
  \draw [color=Orange, -stealth] (\multiplicationFactor*12, \multiplicationFactor*6) -- (\multiplicationFactor*6.1, \multiplicationFactor*13.1);
  \node[color=Orange, right] at (\multiplicationFactor*7.5, \multiplicationFactor*7.5) {$\frac{1}{4}$};
\end{tikzpicture}
  \caption{Prolongation operator in the two-dimensional case with weights.
    The blue arrow represents the operator if $\indexCellFine{i}$ is even and $\indexCellFine{j}$ is even. 
    The red arrows represent the operator if $\indexCellFine{i}$ is odd and $\indexCellFine{j}$ is even.
    The green arrows represent the operator if $\indexCellFine{i}$ is even and $\indexCellFine{j}$ is odd.
  The orange arrows represent the operator if $\indexCellFine{i}$ is odd and $\indexCellFine{j}$ is odd.}
  \label{fig:prolongation2D}
\end{figure}
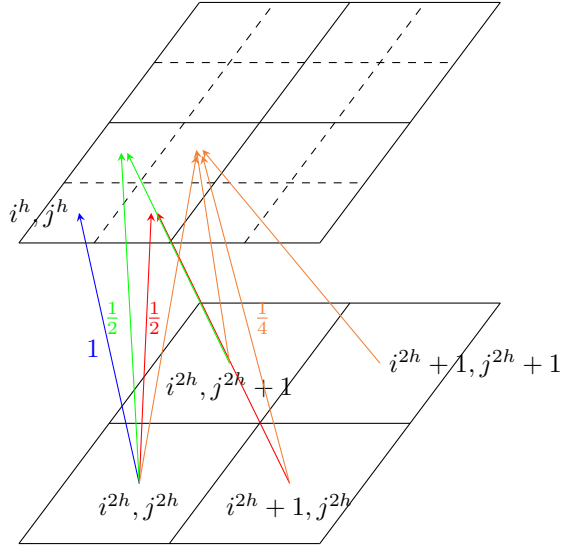

We present now the prolongation operator in the two-dimensional case (see \cref{fig:prolongation2D}).
A cell $\indexCellFine{i}, \indexCellFine{j}$ in the fine mesh $\poissonMeshFine$ can be mapped directly into the cell $\indexCellCoarse{i}, \indexCellCoarse{j}$ in the coarse mesh $\poissonMeshCoarse$, with $\indexCellFine{i} = 2 \indexCellCoarse{i}$ and $\indexCellFine{j} = 2 \indexCellCoarse{j}$.
The value in a cell of the fine mesh is obtained by the interpolation in one direction, followed by the interpolation in the other direction.
The interpolation in the $x$-direction gives
\begin{equation}
    \left( \prolongationOperator{\spaceStep}{2\spaceStep}{\coarseVector{\jacobiUnknown}} \right)_{2 \indexCellCoarse{i}, 2 \indexCellCoarse{j}} = \coarseVector{\jacobiUnknown}_{\indexCellCoarse{i}, \indexCellCoarse{j}} \text{ and } \left( \prolongationOperator{\spaceStep}{2\spaceStep}{\coarseVector{\jacobiUnknown}} \right)_{2 \indexCellCoarse{i}+1, 2 \indexCellCoarse{j}} = \frac{1}{2} \left( \coarseVector{\jacobiUnknown}_{\indexCellCoarse{i}, \indexCellCoarse{j}} + \coarseVector{\jacobiUnknown}_{\indexCellCoarse{i}+1, \indexCellCoarse{j}} \right).
    \label{eq:prolongation_1d_x}
\end{equation}
The interpolation in the $y$-direction preserves these values and gives
\begin{equation}
    \left( \prolongationOperator{\spaceStep}{2\spaceStep}{\coarseVector{\jacobiUnknown}} \right)_{2 \indexCellCoarse{i}, 2 \indexCellCoarse{j}+1} = \frac{1}{2} \left( \coarseVector{\jacobiUnknown}_{\indexCellCoarse{i}, \indexCellCoarse{j}} + \coarseVector{\jacobiUnknown}_{\indexCellCoarse{i}, \indexCellCoarse{j}+1} \right).
    \label{eq:prolongation_1d_y}
\end{equation}
Finally, the last value is given by 
\begin{equation}
    \left( \prolongationOperator{\spaceStep}{2\spaceStep}{\coarseVector{\jacobiUnknown}} \right)_{2 \indexCellCoarse{i}+1, 2 \indexCellCoarse{j}+1} = \frac{1}{4} \left( \coarseVector{\jacobiUnknown}_{\indexCellCoarse{i}, \indexCellCoarse{j}} + \coarseVector{\jacobiUnknown}_{\indexCellCoarse{i+1}, \indexCellCoarse{j}} + \coarseVector{\jacobiUnknown}_{\indexCellCoarse{i}, \indexCellCoarse{j}+1} + \coarseVector{\jacobiUnknown}_{\indexCellCoarse{i}+1, \indexCellCoarse{j}+1} \right) 
    \label{eq:prolongation_1d_average}
\end{equation}
Using \cref{eq:prolongation_1d_x,eq:prolongation_1d_y,eq:prolongation_1d_average} leads to writing the prolongation operator as
\begin{equation}
    \left( \prolongationOperator{\spaceStep}{2\spaceStep}{\coarseVector{\jacobiUnknown}} \right)_{\indexCellFine{i}, \indexCellFine{j}} = 
    \left\{
        \begin{aligned}
            &\coarseVector{\jacobiUnknown}_{\indexCellCoarse{i}, \indexCellCoarse{j}} \text{ if } \indexCellFine{i} \text{ is even and } \indexCellFine{j} \text{ is even},\\
            &\frac{1}{2} \left( \coarseVector{\jacobiUnknown}_{\indexCellCoarse{i}, \indexCellCoarse{j}} + \coarseVector{\jacobiUnknown}_{\indexCellCoarse{i}+1, \indexCellCoarse{j}} \right) \text{ if } \indexCellFine{i} \text{ is odd and } \indexCellFine{j} \text{ is even},\\
            &\frac{1}{2} \left( \coarseVector{\jacobiUnknown}_{\indexCellCoarse{i}, \indexCellCoarse{j}} + \coarseVector{\jacobiUnknown}_{\indexCellCoarse{i}, \indexCellCoarse{j+1}} \right) \text{ if } \indexCellFine{i} \text{ is even and } \indexCellFine{j} \text{ is odd},\\
            &\frac{1}{4} \left( \coarseVector{\jacobiUnknown}_{\indexCellCoarse{i}, \indexCellCoarse{j}} + \coarseVector{\jacobiUnknown}_{\indexCellCoarse{i+1}, \indexCellCoarse{j}} + \coarseVector{\jacobiUnknown}_{\indexCellCoarse{i}, \indexCellCoarse{j}+1} + \coarseVector{\jacobiUnknown}_{\indexCellCoarse{i}+1, \indexCellCoarse{j}+1} \right) \\
            &\quad \text{ if } \indexCellCoarse{i} \text{ is odd and } \indexCellCoarse{j} \text{ is odd}.
        \end{aligned}
    \right.
    \label{eq:prolongationOperator}
\end{equation}

Let us now consider the restriction operator.
We first consider the one-dimensional case, where the mesh is given by \cref{fig:prolongation1D}.
The prolongation operator can be made explicit, with appropriate boundary conditions:
\begin{equation}
    \fineVector{\jacobiUnknown} = \left( \prolongationOperator{\spaceStep}{2\spaceStep}{\coarseVector{\jacobiUnknown}} \right) = 
    \begin{pmatrix}
        \fineVector{\jacobiUnknown}_1\\\fineVector{\jacobiUnknown}_2\\\fineVector{\jacobiUnknown}_3\\\fineVector{\jacobiUnknown}_4\\\fineVector{\jacobiUnknown}_5\\\fineVector{\jacobiUnknown}_6\\\fineVector{\jacobiUnknown}_7
    \end{pmatrix}
    = 
    \begin{pmatrix}
        \fineVector{\jacobiUnknown}_{2 \times 0 + 1}\\ \fineVector{\jacobiUnknown}_{2 \times 1}\\ \fineVector{\jacobiUnknown}_{2 \times 1 + 1}\\ \fineVector{\jacobiUnknown}_{2 \times 2}\\ \fineVector{\jacobiUnknown}_{2 \times 2 + 1}\\ \fineVector{\jacobiUnknown}_{2 \times 3}\\ \fineVector{\jacobiUnknown}_{2 \times 3 + 1}
    \end{pmatrix}
    = \frac{1}{2}
    \begin{pmatrix}
        1 & 0 & 0\\
        2 & 0 & 0\\
        1 & 1 & 0\\
        0 & 2 & 0\\
        0 & 1 & 1\\
        0 & 0 & 2\\
        0 & 0 & 1\\
    \end{pmatrix}
    \begin{pmatrix}
        \coarseVector{\jacobiUnknown}_1\\ \coarseVector{\jacobiUnknown}_2\\ \coarseVector{\jacobiUnknown}_3
    \end{pmatrix}
    .
\end{equation}
Even though the prolongation operator is written as an operator, it is linear and can therefore be seen as a matrix, $\prolongationOperatorNoArg{\spaceStep}{2\spaceStep}$.
Apart from the first and the last rows that handle boundary conditions, the sum of the coefficients in a row is \num{1}.  
Values on the fine mesh are convex combinations of values on the coarse mesh.

We take the restriction operator as the transpose of the prolongation operator: $\restrictionOperatorNoArg{\spaceStep}{2\spaceStep} = \frac{1}{2^d} \transpose{\left( \prolongationOperatorNoArg{\spaceStep}{2\spaceStep} \right)}$, where $d$ is the dimension of the problem ($d \in \left\{ 1,2,3 \right\}$).
Thanks to the factor $\frac{1}{2^d}$, values on the coarse mesh are also convex combinations of values on the fine mesh.
In the one-dimensional case given by \cref{fig:prolongation1D}, the restriction operator writes
\begin{equation}
    \coarseVector{\jacobiUnknown} = \left( \restrictionOperator{\spaceStep}{2\spaceStep}{\coarseVector{\jacobiUnknown}} \right) = 
    \begin{pmatrix}
        \coarseVector{\jacobiUnknown}_1\\ \coarseVector{\jacobiUnknown}_2\\ \coarseVector{\jacobiUnknown}_3
    \end{pmatrix}
    = \frac{1}{4}
    \begin{pmatrix}
        1 & 2 & 1 & 0 & 0 & 0 & 0\\
        0 & 0 & 1 & 2 & 1 & 0 & 0\\
        0 & 0 & 0 & 0 & 1 & 2 & 1\\
    \end{pmatrix}
    \begin{pmatrix}
        \fineVector{\jacobiUnknown}_1\\\fineVector{\jacobiUnknown}_2\\\fineVector{\jacobiUnknown}_3\\\fineVector{\jacobiUnknown}_4\\\fineVector{\jacobiUnknown}_5\\\fineVector{\jacobiUnknown}_6\\\fineVector{\jacobiUnknown}_7
    \end{pmatrix}
    .
\end{equation}
In the general one-dimensional case, we have
\begin{equation}
    \begin{aligned}
        \left( \restrictionOperator{\spaceStep}{2\spaceStep}{\fineVector{\jacobiUnknown}} \right)_{\indexCellCoarse{i}} &=  
        \frac{1}{4} \fineVector{\jacobiUnknown}_{\indexCellFine{i}-1}
        + \frac{1}{2} \fineVector{\jacobiUnknown}_{\indexCellFine{i}  }
        + \frac{1}{4} \fineVector{\jacobiUnknown}_{\indexCellFine{i}+1}.
    \end{aligned}
\end{equation}
Numerical tests in \cref{sect:gmg:resultatsNum} are done in the two-dimensional case, we also make $\restrictionOperatorNoArg{\spaceStep}{2\spaceStep}$ explicit with $d = 2$:
\begin{equation}
    \begin{aligned}
        \left( \restrictionOperator{\spaceStep}{2\spaceStep}{\fineVector{\jacobiUnknown}} \right)_{\indexCellCoarse{i}, \indexCellCoarse{j}} &=  
        \frac{1}{16} \fineVector{\jacobiUnknown}_{\indexCellFine{i}-1, \indexCellFine{j}-1}
        + \frac{1}{8 } \fineVector{\jacobiUnknown}_{\indexCellFine{i}-1, \indexCellFine{j}  }
        + \frac{1}{16} \fineVector{\jacobiUnknown}_{\indexCellFine{i}-1, \indexCellFine{j}+1}\\
        & + \frac{1}{8 } \fineVector{\jacobiUnknown}_{\indexCellFine{i}  , \indexCellFine{j}-1}
        + \frac{1}{4 } \fineVector{\jacobiUnknown}_{\indexCellFine{i}  , \indexCellFine{j}  }
        + \frac{1}{8 } \fineVector{\jacobiUnknown}_{\indexCellFine{i}  , \indexCellFine{j}+1}\\
        & + \frac{1}{16} \fineVector{\jacobiUnknown}_{\indexCellFine{i}+1, \indexCellFine{j}-1}
        + \frac{1}{8 } \fineVector{\jacobiUnknown}_{\indexCellFine{i}+1, \indexCellFine{j}  }
        + \frac{1}{16} \fineVector{\jacobiUnknown}_{\indexCellFine{i}+1, \indexCellFine{j}+1}.
    \end{aligned}
    \label{eq:restrictionOperator}
\end{equation}

Now that we have made explicit the restriction and prolongation operators used with two grids, we write the two-grid full approximation scheme.

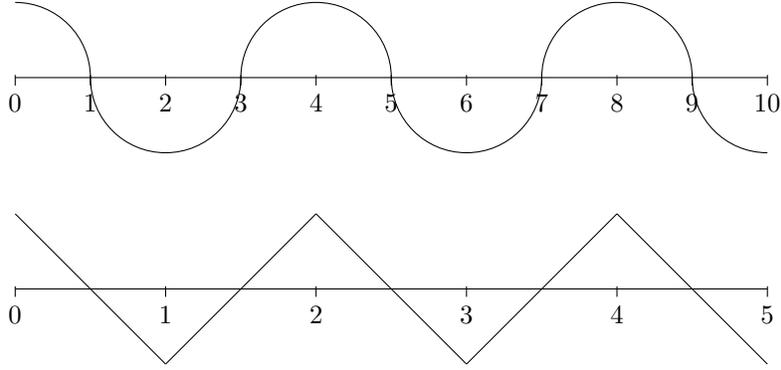
\begin{figure}
    \centering
    \begin{tikzpicture}
  \begin{scope}[yshift=80] 
    \draw (0,0) -- (10,0);
    \foreach \x in {0,...,10} \draw (\x,1pt) -- (\x,-3pt) node[anchor=north] {\x};
    \draw (1,0) arc (0:90:1);
    \draw (3,0) arc (0:-180:1);
    \draw (5,0) arc (0:180:1);
    \draw (7,0) arc (0:-180:1);
    \draw (9,0) arc (0:180:1);
    \draw (9,0) arc (-180:-90:1);
  \end{scope}
  \begin{scope}
    \draw (0,0) -- (10,0);
    \foreach \x in {0,...,5} \draw (2*\x,1pt) -- (2*\x,-3pt) node[anchor=north] {\x};
    \draw (0,1) -- (2,-1);
    \draw (2,-1) -- (4,1);
    \draw (4,1) -- (6,-1);
    \draw (6,-1) -- (8,1);
    \draw (8,1) -- (10,-1);
  \end{scope}
\end{tikzpicture}
    \caption[Restriction of a smooth function]{Restriction of a smooth function.}
    \label{fig:gmg:linear:smooth}
\end{figure}

One can show that the presence of smooth modes in the error between the approximate and exact solution can cause the relaxation to stall.
In the linear case, the proof relies on the study of the eigenvalues and eigenvectors of $\jacobiOperator$, see e.g., \cite{brandt2011} where the Poisson equation is solved.
As shown by \cref{fig:gmg:linear:smooth}, a smooth function on the fine mesh will be less smooth on the coarse mesh. 
An iterative method such as Jacobi will converge faster and will smooth the error on the coarse mesh.
This gives the idea of restricting and prolonging the error between the different meshes.
The function plotted in \cref{fig:gmg:linear:smooth} is a generic function to illustrate the interest of the restriction.
It does not relate with our problem.
If $\vector{\jacobiUnknown}$ is an approximate solution of \cref{eq:gmg:nonLinearSyst} and $\vector{\jacobiIntermediate}$ is the exact solution, let us define the error as $\vector{\jacobiError} = \vector{\jacobiIntermediate} - \vector{\jacobiUnknown}$.
By defining the residual $\vector{\jacobiResidual} = \vector{\jacobiRhs} - \jacobiOperator \vector{\jacobiUnknown}$ (\cref{eq:defResidual}), the error verifies the following equation:
\begin{equation}
    \begin{aligned}
        \vector{\jacobiResidual} &= \vector{\jacobiRhs} - \jacobiOperator \left( \vector{\jacobiUnknown} \right)\\
        &= \jacobiOperator \left( \vector{\jacobiIntermediate} \right) - \jacobiOperator \left( \vector{\jacobiUnknown} \right).
    \end{aligned}
    \label{eq:residualEquaNonLinear}
\end{equation}
The error can be computed by solving the residual equation \cref{eq:residualEquaNonLinear} on the coarse mesh:
\begin{equation}
    \coarseOperator \left( \coarseVector{\jacobiIntermediate} \right) = \coarseOperator \left( \coarseVector{\jacobiUnknown} \right) + \coarseVector{\jacobiResidual},
    \label{eq:equaNonLinearCoarseGrid}
\end{equation}
whose right-hand side is known.
$\coarseVector{\jacobiResidual} = \restrictionOperator{\spaceStep}{2\spaceStep}{\fineVector{\jacobiResidual}} = \restrictionOperator{\spaceStep}{2\spaceStep}{\fineVector{\jacobiRhs} - \fineOperator\left( \fineVector{\jacobiUnknown} \right)}$ is an approximation of the residual on the coarse mesh.
Likewise, $\coarseVector{\jacobiUnknown} = \restrictionOperator{\spaceStep}{2\spaceStep}{\fineVector{\jacobiUnknown}}$.
The new approximate solution on the fine mesh can be computed:
\begin{equation}
    \fineVector{\jacobiUnknown} + \prolongationOperator{\spaceStep}{2\spaceStep}{\underbrace{\coarseVector{\jacobiIntermediate} - \coarseVector{\jacobiUnknown}}_{\coarseVector{\jacobiError}}}.
    \label{eq:correctionFineGrid}
\end{equation}
The solution of the problem solved on the coarse mesh is the full approximation $\coarseVector{\jacobiIntermediate} = \coarseVector{\jacobiUnknown} + \coarseVector{\jacobiError}$, and not the error $\coarseVector{\jacobiError}$.
This leads to the full approximation scheme \cref{algo:fas} (\cite{briggs2000}).

One could try to solve $\coarseOperator \left( \coarseVector{\jacobiUnknown} \right) = \coarseVector{\jacobiRhs}$ on the fine grid instead of solving \cref{eq:equaNonLinearCoarseGrid}.
However, numerical experiments show that solving \cref{eq:equaNonLinearCoarseGrid} and correcting the solution with \cref{eq:correctionFineGrid} is required to achieve convergence on the fine grid.

\begin{algorithm}
    \caption{Two-grid full approximation scheme}
    \begin{algorithmic}[1]
        \State Pre-smoother: relax $\preSmootherNbIter$ times $\fineOperator \left( \fineVector{\jacobiIntermediate} \right) = \fineVector{\jacobiRhs}$ with initial guess $\fineVector{\jacobiUnknown}$ on the fine mesh $\poissonMeshFine$
        \State Restrict the residual $\coarseVector{\jacobiResidual} = \restrictionOperator{\spaceStep}{2\spaceStep}{\fineVector{\jacobiRhs} - \fineOperator \left( \fineVector{\jacobiIntermediate} \right)}$ and the current solution $\coarseVector{\jacobiUnknown} = \restrictionOperator{\spaceStep}{2\spaceStep}{\fineVector{\jacobiIntermediate}}$
        \State Solve $\coarseOperator \left( \coarseVector{\jacobiIntermediate} \right) = \coarseVector{\jacobiResidual} + \coarseOperator \left( \coarseVector{\jacobiUnknown} \right)$ with initial guess $\coarseVector{\jacobiUnknown}$ on the coarse mesh $\poissonMeshCoarse$
        \State Correct the approximation on the fine mesh: $\fineVector{\bar{\jacobiUnknown}} = \fineVector{\jacobiUnknown} + \prolongationOperator{\spaceStep}{2\spaceStep}{\coarseVector{\jacobiIntermediate} - \coarseVector{\jacobiUnknown}}$
        \State Post-smoother: relax $\postSmootherNbIter$ times $\fineOperator \left( \fineVector{\jacobiIntermediate} \right) = \fineVector{\jacobiRhs}$ with initial guess $\fineVector{\bar{\jacobiUnknown}}$ on the fine mesh $\poissonMeshFine$
    \end{algorithmic}
    \label{algo:fas}
\end{algorithm}

If $\jacobiOperator$ were linear, we would have $\vector{\jacobiResidual} = \jacobiOperator \left( \vector{\jacobiIntermediate} \right) - \jacobiOperator \left( \vector{\jacobiUnknown} \right) = \jacobiOperator \left( \vector{\jacobiIntermediate} - \vector{\jacobiUnknown} \right) = \jacobiOperator \left( \vector{\jacobiError} \right)$.
Instead of solving \cref{eq:equaNonLinearCoarseGrid} on the coarse mesh, we would solve
\begin{equation}
  \coarseOperator \left( \coarseVector{\jacobiError} \right) = \coarseVector{\jacobiResidual},
\end{equation}
and \cref{eq:correctionFineGrid} would be replaced by
\begin{equation}
  \fineVector{\jacobiUnknown} + \prolongationOperator{\spaceStep}{2\spaceStep}{\coarseVector{\jacobiError}}.
\end{equation}
Because $\jacobiOperator$ is nonlinear, we have to use the FAS algorithm.

This process can now be applied recursively until a mesh coarse enough is reached.
This leads to the so-called V-cycle algorithm.
Solving the nonlinear system at a given level is less costly than the previous level because the number of unknowns is reduced.
However, numerical experiments have shown that the best performances are obtained by performing a series of V-cycles.
In the linear case, a direct solver could be used at the coarsest level.
In the case of radiative transfer, we will use the Jacobi method as a coarse grid solver.

\subsection{Application to the HLL solver for the M\textsubscript{1} model} \label{sect:gmg:m1}

Let us now apply \cref{algo:fas} to solve \cref{eq:gmg:jacobi:hll}.
The operator $\jacobiOperator$ is the one described in \cref{sect:gmg:jacobi}, and we use \cref{algo:jacobi} as smoother and coarse grid solver.
The system to be solved on the coarse mesh is 
\begin{equation}
    \coarseOperator \left( \coarseVector{\jacobiIntermediate} \right) = \coarseOperator \left( \coarseVector{\jacobiUnknown} \right) + \coarseVector{\jacobiResidual}.
    \label{eq:gmg:pseudoTime:initialRhs}
\end{equation}
As shown in \cref{sect:jacobi:etatsAdmissibles}, if the initial guess for $\coarseVector{\jacobiIntermediate}$ and the right-hand side $\coarseOperator \left( \coarseVector{\jacobiUnknown} \right) + \coarseVector{\jacobiResidual}$ are admissible, then the solution obtained with \cref{algo:jacobi} is also admissible.
However, numerical experiments have shown that, in general, $\coarseOperator \left( \coarseVector{\jacobiUnknown} \right) + \coarseVector{\jacobiResidual}$ is not admissible, which leads to a non-admissible solution.
If the right-hand side of \cref{eq:gmg:pseudoTime:initialRhs} is not admissible, the Jacobi method will theoretically converge to a unique solution, this solution might not be admissible.

To tackle this issue, we follow \cite{kifonidis2012} and we introduce a pseudo-time $\pseudoTime$. 
Instead of solving \cref{eq:gmg:pseudoTime:initialRhs}, we look for the steady state in pseudo-time of the following equation:
\begin{equation}
    \frac{\dif \coarseVector{\jacobiIntermediate}}{\dif \pseudoTime} + \coarseOperator \left( \coarseVector{\jacobiIntermediate} \right) = \coarseOperator \left( \coarseVector{\jacobiUnknown} \right) + \coarseVector{\jacobiResidual}.
    \label{eq:gmg:pseudoTime:pseudoTimeEqua}
\end{equation}
When the steady state is reached, $\frac{\dif \coarseVector{\jacobiIntermediate}}{\dif \pseudoTime} = 0$, and we recover \cref{eq:gmg:pseudoTime:initialRhs}. 
Let us notice that the pseudo-time $\pseudoTime$ is completely independent of the physical time step $\Delta t$.

\Cref{eq:gmg:pseudoTime:pseudoTimeEqua} is a (nonlinear) system of ordinary differential equations in the variable $\pseudoTime$.
We use a notation similar to the physical time for the discretization in pseudo-time.
$\Delta \pseudoTime$ is the interval between the current pseudo-time $\pseudoTime^{\pseudoTimeIndex}$ and $\pseudoTime^{\pseudoTimeIndex+1}$.
We choose $\Delta \pseudoTime$ such that $\coarseVector{\jacobiIntermediate}$ is admissible.
Using the definition of the residual, the right-hand side of \cref{eq:gmg:pseudoTime:pseudoTimeEqua} becomes $\coarseOperator \left( \coarseVector{\jacobiUnknown} \right) + \coarseVector{\jacobiRhs} - \restrictionOperator{\spaceStep}{2\spaceStep}{\fineOperator \left( \fineVector{\jacobiUnknown} \right)}$.
We want to solve \cref{eq:gmg:pseudoTime:pseudoTimeEqua} with a stable scheme, for all $\Delta \pseudoTime$.
Therefore, we use a splitting strategy.
For the scheme to be stable, the left-hand side has to be taken implicitly.
This leads to solving \cref{eq:gmg:pseudoTime:pseudoTimeEqua} as 
\begin{equation}
    \begin{aligned}
        \frac{\widetilde{\coarseVector{\jacobiIntermediate}} - \left( \coarseVector{\jacobiIntermediate} \right)^{\pseudoTimeIndex}}{\Delta \pseudoTime} &= \coarseOperator \left( \coarseVector{\jacobiUnknown} \right) + \coarseVector{\jacobiRhs} - \restrictionOperator{\spaceStep}{2\spaceStep}{\fineOperator \left( \fineVector{\jacobiUnknown} \right)}\\
        \frac{\left( \coarseVector{\jacobiIntermediate} \right)^{\pseudoTimeIndex+1} - \widetilde{\coarseVector{\jacobiIntermediate}}}{\Delta \pseudoTime} &+ \coarseOperator \left( \left( \coarseVector{\jacobiIntermediate} \right)^{\pseudoTimeIndex+1} \right) = 0.\\
    \end{aligned}
    \label{eq:gmg:pseudoTime:splitting}
\end{equation}

\begin{propostion}
  $\left( \coarseVector{\jacobiIntermediate} \right)^{\pseudoTimeIndex+1}$ given by \cref{eq:gmg:pseudoTime:splitting} is admissible.
\end{propostion}

\begin{proof}
  The first equation in \cref{eq:gmg:pseudoTime:splitting} is explicit in pseudo-time.
  We can always choose a value for $\Delta \pseudoTime$ such that $\widetilde{\coarseVector{\jacobiIntermediate}}$ is admissible.
  Let us notice that the right-hand side is fixed and the left-hand side is local to a cell. 

  The second equation in \cref{eq:gmg:pseudoTime:splitting} is implicit in pseudo-time.
  With arguments similar to \cref{sect:jacobi:etatsAdmissibles}, $\left( \coarseVector{\jacobiIntermediate} \right)^{\pseudoTimeIndex+1}$ is admissible as soon as $\widetilde{\vector{\jacobiIntermediate}^{2h}}$ is admissible.
\end{proof}

Unfortunately, choosing $\Delta \pseudoTime$ such that $\widetilde{\coarseVector{\jacobiIntermediate}}$ is admissible can result in small values for $\Delta \pseudoTime$.
To reduce the computational cost, we slightly change the algorithm.
We do not use the same pseudo-time step in all cells and for the implicit step.
We choose a pseudo-time step $\dTauIm$ for the implicit step, and we use a sub-cycle algorithm with a time step local to each cell until $\dTauIm$ is reached.
This results in \cref{algo:pseudoTimeTwoGrids}.
We use this algorithm as a smoother and coarse solver in \cref{algo:fas}.

\begin{algorithm}
    \caption{Resolution at coarse level using pseudo-time}
    \begin{algorithmic}[1]
        \State Choose $\dTauIm$
        \State $\pseudoTimeIndex = 0$
        \While{the steady state is not reached}
        \For{each cell $\indexCellCoarse{i}$ in the coarse grid}
        \State $\tauEx_{\indexCellCoarse{i}} = 0$
        \State $\pseudoTimeCounter = 0$
        \State $\left( \widetilde{\coarseVector{\jacobiIntermediate}} \right)^{(0)}_{\indexCellCoarse{i}} = \left( \coarseVector{\jacobiIntermediate} \right)^{\pseudoTimeIndex}_{\indexCellCoarse{i}}$
        \While{$\tauEx_{\indexCellCoarse{i}} < \dTauIm$}
        \State $\left( \widetilde{\coarseVector{\jacobiIntermediate}} \right)^{(\pseudoTimeCounter+1)}_{\indexCellCoarse{i}} = \left( \widetilde{\coarseVector{\jacobiIntermediate}} \right)^{(\pseudoTimeCounter)}_{\indexCellCoarse{i}} + \dTauEx_{\indexCellCoarse{i}} \left( \left( \coarseOperator \left( \coarseVector{\jacobiUnknown} \right) \right)_{\indexCellCoarse{i}} + \coarseVector{\jacobiRhs}_{\indexCellCoarse{i}} - \left( \restrictionOperator{\spaceStep}{2\spaceStep}{\fineOperator \left( \fineVector{\jacobiUnknown} \right)} \right)_{\indexCellCoarse{i}} \right)$ with $\dTauEx_{\indexCellCoarse{i}}$ such that $\left( \widetilde{\coarseVector{\jacobiIntermediate}}\right)^{(\pseudoTimeCounter+1)}_{\indexCellCoarse{i}}$ is admissible
        \State $\pseudoTimeCounter \gets \pseudoTimeCounter+1$
        \State $\tauEx_{\indexCellCoarse{i}} \gets \tauEx_{\indexCellCoarse{i}} + \dTauEx_{\indexCellCoarse{i}}$
        \EndWhile
        \State $\widetilde{\coarseVector{\jacobiIntermediate}}_{\indexCellCoarse{i}} = \left( \widetilde{\coarseVector{\jacobiIntermediate}}\right)^{(\pseudoTimeCounter)}_{\indexCellCoarse{i}}$
        \EndFor
        \State Solve $\left( \coarseVector{\jacobiIntermediate} \right)^{\pseudoTimeIndex+1} + \dTauIm \coarseOperator \left( \left( \coarseVector{\jacobiIntermediate} \right)^{\pseudoTimeIndex+1} \right) = \widetilde{\coarseVector{\jacobiIntermediate}}$ with nonlinear Jacobi method (\cref{algo:jacobi})
        \State $\pseudoTimeIndex \gets \pseudoTimeIndex + 1$
        \EndWhile
    \end{algorithmic}
    \label{algo:pseudoTimeTwoGrids}
\end{algorithm}

The next question to be solved is how to choose $\dTauIm$.
Numerical experiments have shown that a large value for $\dTauIm$ allows reducing the computational cost per V-cycle, but can result in a very low decrease of the norm of the residual.
To avoid such a result, we use an adaptive pseudo-time step.
When the norm of the residual decreases fast enough, $\dTauIm$ is increased.
On the contrary, when the norm of the residual decreases slowly, $\dTauIm$ is decreased.

$\dTauEx$ is chosen small enough so that each component of $\widetilde{\coarseVector{\jacobiIntermediate}}$ is admissible.
Using the first equation of \cref{eq:gmg:pseudoTime:splitting}, one can obtain an explicit formula for $\dTauEx$.
However, this can result is very small values of $\dTauEx$.
Because we use this algorithm only for the coarse levels, numerical experiments show that the performing only a few iterations is enough to improve the convergence of the method on the finest level (see \cref{sect:gmg:resultatsNum}).

Furthermore, with this choice of $\dTauEx$ and $\dTauIm$, the Jacobi method used to solve the second equation of \cref{eq:gmg:pseudoTime:splitting}  has always converged in numerical experiments.

Let us emphasize that we do not introduce the pseudo-time $\pseudoTime$ to improve the performances of the method, but preserve the admissible states $\radiativeEnergy > 0$ and $\reducedFlux \le 1$.
Numerical values for $\dTauEx$ and $\dTauIm$ are tuned to have reasonable perforances despite the introduction of the new iterative process.

The update of the solution with the correction on the fine grid (line \num{4} in \cref{algo:fas})
\begin{equation}
    \overline{\fineVector{\jacobiUnknown}} = \fineVector{\jacobiUnknown} + \prolongationOperator{\spaceStep}{2\spaceStep}{\coarseVector{\jacobiIntermediate} - \coarseVector{\jacobiUnknown}}
    \label{eq:gmg:pseudoTime:prolongation}
\end{equation}
can also result in non-admissible states.
Let us introduce another pseudo-time $\pseudoTimeProlongation$ to compute the steady state of
\begin{equation}
    \frac{\dif \overline{\fineVector{\jacobiUnknown}}}{\dif \pseudoTimeProlongation} + \overline{\fineVector{\jacobiUnknown}} = \fineVector{\jacobiUnknown} + \prolongationOperator{\spaceStep}{2\spaceStep}{\coarseVector{\jacobiIntermediate} - \coarseVector{\jacobiUnknown}},
    \label{eq:gmg:pseudoTime:steadyStateProlongation}
\end{equation}
instead of using \cref{eq:gmg:pseudoTime:prolongation}.
We write $\Delta \pseudoTimeProlongation$ the interval between the current pseudo-time $\pseudoTimeProlongation^{\overline{\pseudoTimeIndex}}$ and $\pseudoTimeProlongation^{\overline{\pseudoTimeIndex}+1}$.
We can now discretize \cref{eq:gmg:pseudoTime:steadyStateProlongation} with a time-implicit solver:
\begin{equation}
    \frac{\left( \overline{\fineVector{\jacobiUnknown}} \right)^{\overline{\pseudoTimeIndex}+1} - \left( \overline{\fineVector{\jacobiUnknown}} \right)^{\overline{\pseudoTimeIndex}}}{\Delta \pseudoTimeProlongation} + \left( \overline{\fineVector{\jacobiUnknown}} \right)^{\overline{\pseudoTimeIndex}+1} = \fineVector{\jacobiUnknown} + \prolongationOperator{\spaceStep}{2\spaceStep}{\coarseVector{\jacobiIntermediate} - \coarseVector{\jacobiUnknown}}. 
\end{equation}
We choose $\Delta \pseudoTimeProlongation$ such that
\begin{equation}
    \left( \overline{\fineVector{\jacobiUnknown}} \right)^{\overline{\pseudoTimeIndex}+1} = \frac{\left( \overline{\fineVector{\jacobiUnknown}} \right)^{\overline{\pseudoTimeIndex}} + \Delta \pseudoTimeProlongation \left( \fineVector{\jacobiUnknown} + \prolongationOperator{\spaceStep}{2\spaceStep}{\coarseVector{\jacobiIntermediate} - \coarseVector{\jacobiUnknown}} \right)}{1 + \Delta \pseudoTimeProlongation}
    \label{eq:gmg:pseudoTime:prolongatedSol}
\end{equation}
is admissible.
This process can be applied locally, only in the cells where it is needed.

The whole algorithm, including these modifications, is presented in \labelcref{appendix:gmg:overallAlgo}.

In the next section, we present numerical results to show the gain in computational time.

\section{Numerical results} \label{sect:gmg:resultatsNum}

We perform some verification tests to validate the performances of the algorithms presented in \cref{sect:gmg:jacobi,sect:gmg:multigrille}.
We compare our results with those obtained with a time-explicit HLL solver and a time-implicit HLL solver where a Newton-Raphson method is used to solve \cref{eq:gmg:jacobi:hll}.
We will compare the performances with different time steps.
Furthermore, we write $\Delta t = \cfl \frac{\spaceStep}{\speedOfLight}$.
For the time-explicit solver, one should have $\cfl \le 1$ to respect the CFL condition.

We always use the same parameters for the geometric multigrid method. At the finest level, the number of iterations for the pre- and post-smoothers is $\smootherNbIter{0} = 3$. 
When using iterations in pseudo-time, the number of iterations in pseudo-time $\pseudoTimeIndex$ is set to \num{3} and the number of iterations for the smoothers is $\smootherNbIter{\indexLevel} = 1$, for $\indexLevel \neq 0$.
The initial value for $\dTauIm$ is \num{e-3}.
These parameters are chosen because they give reasonable performances in most cases and can easily be used in physical problems.

Computational times shown later are obtained with a sequential implementation.
All tests were run with an Intel Core i7-8650U @ \SI{1.90}{\giga \hertz} processor.

\subsection{Beam} \label{sect:gmg:beam}

\begin{figure}
    \begin{subfigure}[b]{0.5\textwidth}
        \begin{centering}
            \includegraphics[width=\columnwidth]{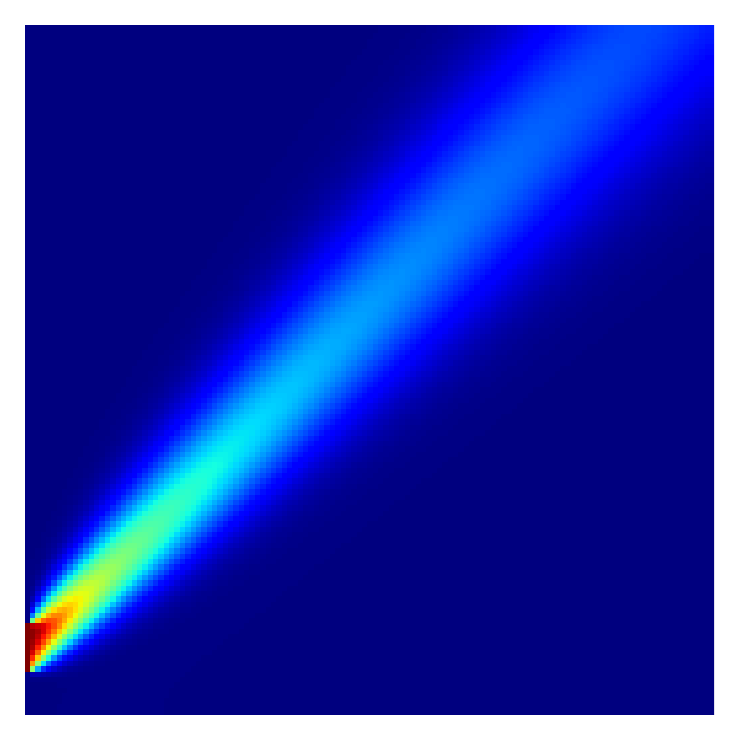}
            \caption{Explicit}
            \label{fig:gmg:beam:explicit}
        \end{centering}
    \end{subfigure}
    \begin{subfigure}[b]{0.5\textwidth}
        \begin{centering}
            \includegraphics[width=\columnwidth]{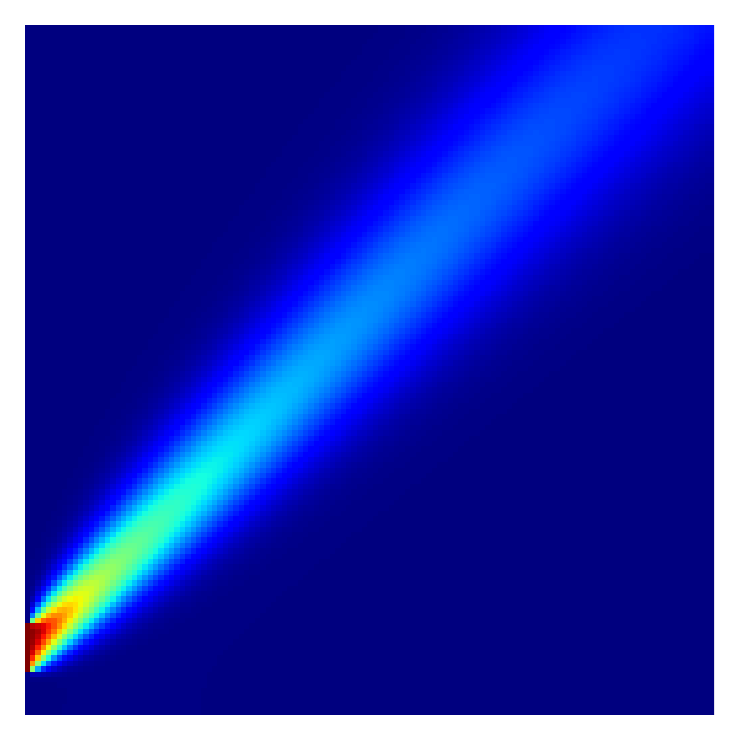}
            \caption{Jacobi method}
            \label{fig:gmg:beam:jacobi}
        \end{centering}
    \end{subfigure}
    \caption{Snapshots of radiative energy at steady state with the explicit solver and the Jacobi method.}
    \label{fig:gmg:beam}
\end{figure}

\begin{figure}
    \centering
    \begin{tikzpicture}
  \begin{axis}[
      xlabel = {$x$ (\si{\centi \meter})},
      ylabel = {Radiative energy (arbitrary unit)},
      ymax = {1.01},
      axis lines = left,
    ]
    \addplot[color=black] table {images/beam/cut_hll.txt};
    \addlegendentry{Explicit, $\cfl = \num{0.45}$};
    
    \addplot[color=green] table {images/beam/cut_newton.txt};
    \addlegendentry{Newton-Raphson, $\cfl = \num{0.2}$};
    
    \addplot[color=blue] table {images/beam/cut_level1_2000.txt};
    \addlegendentry{GMG, $\nbLevelMax = 1$, $\cfl = \num{2000}$};

    \addplot[color=red, dotted] table {images/beam/cut_level2_2000.txt};
    \addlegendentry{GMG, $\nbLevelMax \in \{2,3,4\}$, $\cfl = \num{2000}$};

  \end{axis}
\end{tikzpicture}
    \caption[Horizontal cut in beam test using the Jacobi method and the GMG algorithm]{Beam simulation. The figure shows a horizontal cut at the middle height with the explicit solver, the implicit solver using Newton-Raphson method and the geometric multigrid algorithm different values for $\nbLevelMax$.}
    \label{fig:gmg:beam_cut}
\end{figure}
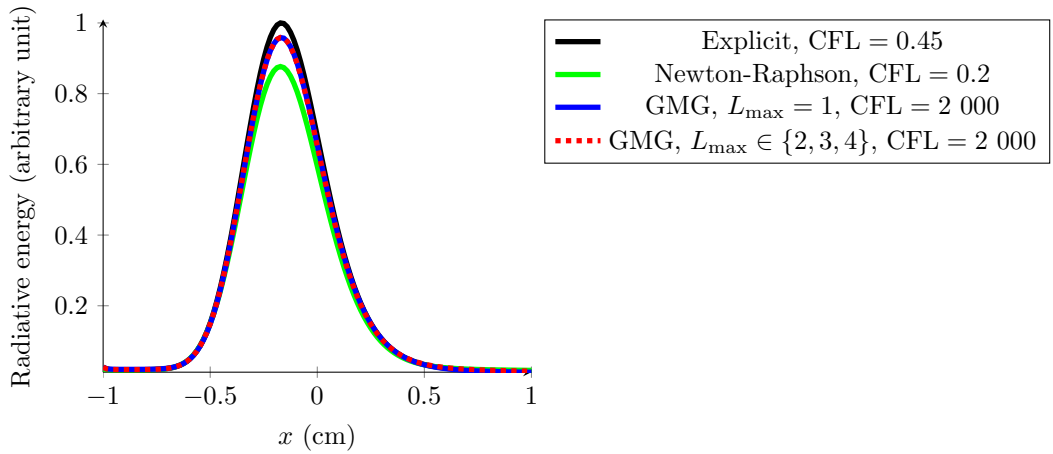

Let us consider the test described in \cref{sect:gmg:interet}.
It is the propagation of a beam in the free-streaming regime.
The domain is first discretized with $129 \times 129$ cells.
\Cref{fig:gmg:beam} shows the radiative energy at steady state.
The solution is obtained with the explicit solver (\cref{fig:gmg:beam:explicit}) and with the Jacobi method (\cref{fig:gmg:beam:jacobi}).
The solution obtained with the geometric multigrid method is similar to \cref{fig:gmg:beam:jacobi}.
\Cref{fig:gmg:beam_cut} shows a horizontal cut at the middle height once the steady state is reached.
We expect the beam to cross the box without dispersion. 
Because we use a scheme of order \num{1}, the solution we obtain is diffused.
We use the solution obtained with the explicit solver as reference to quantify the numerical diffusion introduced by the implicit methods.
Numerical experiments show that the iterations of the Newton-Raphson method do not preserve the admissible states when large time steps are used in the free-streaming regime.
Using the Newton-Raphson method, one has to solve a nonlinear system of the form $\newtonOperator \left( \newtonUnknown \right) = 0$.
At iteration $\newtonIter$, the linear system $\left( \newtonJac{\newtonIter} \right) \newtonError^{(\newtonIter)} = - \newtonOperator \left( \newtonUnknown^{(\newtonIter)} \right)$ is solved.
$\newtonJac{\newtonIter}$ is the Jacobian of $\newtonOperator$.
Then, $\newtonUnknown^{(\newtonIter+1)} = \newtonUnknown^{(\newtonIter)} + \newtonError^{(\newtonIter)}$.
In our case, $\newtonJac{\newtonIter}$ is a large nonsymmetric sparse matrix.
Classical iteration methods such as preconditioned biconjugate gradient stabilized (\cite{vanDerVorst1992}) are not designed to preserve the admissible states.
In the general case, the elements of $\newtonError^{(\newtonIter)}$ are not admissible states, then the elements of $\newtonUnknown^{(\newtonIter+1)}$ are not admissible states either.
The time step has to be chosen small enough to ensure that the elements of $\newtonUnknown^{(\newtonIter+1)}$ are admissible states.
Indeed, we have to use a smaller time step than the one required by the explicit scheme to preserve the admissible states.
We do not expect the M\textsubscript{1} model to perform well in the free-streaming regime with large time steps.
However, we consider this case because the closure relation at the continuous level is chosen such that the model capture the physical behavior in this regime.
Despite reducing the time step,
the solution obtained with the Newton-Raphson method is more diffused than the one obtained with the explicit scheme.
It reaches only $89\%$ of the maximum value of radiative energy.
On the contrary, the solutions obtained with the Jacobi method ($\nbLevelMax = 1$) and with the geometric multigrid algorithm ($\nbLevelMax=2$) reach $96\%$ of the maximum of radiative energy, with a much larger time step ($\cfl = \num{2000}$).

\begin{table*}
  \centering
  \begin{tabular}{cc}     
    \hline
    Scheme      &  Computational time (normalized) \\
    \hline
    Explicit & \num{3} \\
    Newton-Raphson & \num{145} \\
    GMG, $\nbLevelMax = 1$ & \num{1} \\
    GMG, $\nbLevelMax = 2$ & \num{0.85} \\
    GMG, $\nbLevelMax = 3$ & \num{0.79} \\
    GMG, $\nbLevelMax = 4$ & \num{0.38} \\
  \end{tabular}
  \caption[Computation time for the beam problem using the Jacobi method and the GMG algorithm]{Computation time to reach the steady state with the explicit solver, the implicit solver using Newton-Raphson method and the geometric multigrid algorithm different values for $\nbLevelMax$.}
  \label{table:gmg:perfBeam}
\end{table*}

Using the geometric multigrid method should reduce the computational cost.
If we reach a low residual such as $\frac{\norm{\cycleIter{\fineVector{\jacobiResidual}}{\indexNbCycle}}}{\norm{\cycleIter{\fineVector{\jacobiResidual}}{0}}} = \num{e-5}$, numerical experiments have shown that there is no gain in computational cost.
Therefore, we set it to \num{e-2} and we check some properties of the scheme in \cref{sect:gmg:riemannPb}.

\Cref{table:gmg:perfBeam} shows the computational time needed to reach the steady state with different methods: the explicit HLL solver, the implicit HLL solver using the Newton-Raphson method to solve the nonlinear solver, and the geometric multigrid algorithm with different values for $\nbLevelMax$.
The resolution used is now $257 \times 257$ cells.
Using the explicit solver and the Newton-Raphson method, the time step is restricted by the CFL condition, whereas the steady state is reached with only one iteration using the geometric multigrid method.
With $\nbLevelMax = 1$, we recover the Jacobi method, and the steady state is reached three times faster than using the explicit solver.
Increasing the value of $\nbLevelMax$ leads again to a decrease in computational cost.
The time needed to reach the steady state with $\nbLevelMax = 1$ is more than twice the time needed with $\nbLevelMax = 4$.

\begin{table}
  \centering
  \begin{tabular}{cc}
    \hline
    $\nbLevelMax$ & Computational time per V-cycle (normalized)\\
    \hline
    \num{1} & \num{1} \\
    \num{2} & \num{2.9} \\
    \num{3} & \num{3.5} \\
    \num{4} & \num{3.6} \\
  \end{tabular}
  \caption[Computational time per V-cycle for the beam problem]{Computational time per V-cycle with different values of $\nbLevelMax$.}
  \label{tab:gmg:costPerCycle}
\end{table}

\begin{figure}
    \centering
    \begin{tikzpicture}
  \begin{semilogyaxis}[
      xlabel = {Number of V-cycles ($\indexNbCycle$)},
      ylabel = {Residual $\frac{\norm{\cycleIter{\fineVector{\jacobiResidual}}{\indexNbCycle}}}{\norm{\cycleIter{\fineVector{\jacobiResidual}}{0}}}$},
      axis lines = left,
    ]
    \addplot[color=blue] table {images/cvGmgBeam/output_level1.txt};
    \addlegendentry{$\nbLevelMax = 1$};

    \addplot[color=red] table {images/cvGmgBeam/output_level2.txt};
    \addlegendentry{$\nbLevelMax = 2$};

    \addplot[color=green] table {images/cvGmgBeam/output_level3.txt};
    \addlegendentry{$\nbLevelMax = 3$};

    \addplot[color=Orange] table {images/cvGmgBeam/output_level4.txt};
    \addlegendentry{$\nbLevelMax = 4$};

  \end{semilogyaxis}
\end{tikzpicture}
    \caption[Residual as a function of the number of V-cycle for the beam test]{Evolution of the residual as a function of the number of V-cycle, with different values for $\nbLevelMax$ to reach the steady state for the beam problem.}
    \label{fig:gmg:beam:res}
\end{figure}
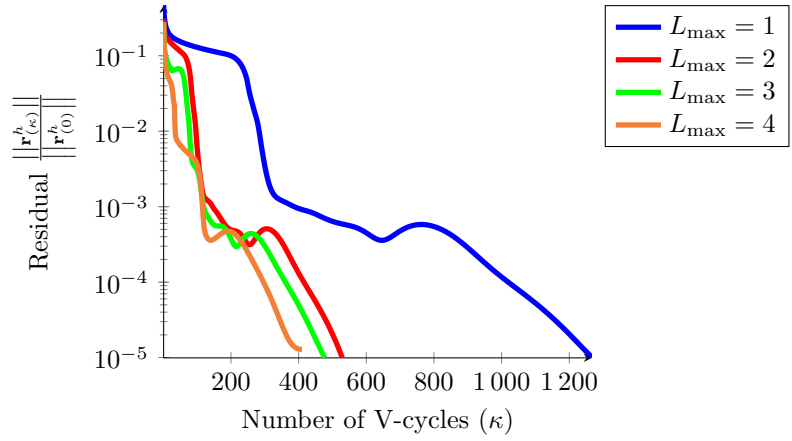

As shown by \cref{tab:gmg:costPerCycle}, when $\nbLevelMax$ increases, the computational time per V-cycle also increases, from \SI{0.48}{\second} with $\nbLevelMax = 1$ to \SI{1.75}{\second} with $\nbLevelMax = 4$.
With $\nbLevelMax > 1$, the computational time per V-cycle is more than three times the computational time per cycle with $\nbLevelMax = 1$.
However, the total computational time (\cref{table:gmg:perfBeam}) decreases when $\nbLevelMax$ increases.
As shown by \cref{fig:gmg:beam:res}, the number of V-cycles to reach the same residual decreases as $\nbLevelMax$ increases, from \num{1267} with $\nbLevelMax = 1$ to \num{400} with $\nbLevelMax = 4$.
With different values for $\nbLevelMax$, all the curves have the same shape, but they are stiffer when $\nbLevelMax$ increases.

\begin{table}
    \centering
    \begin{tabular}{cc}
        \hline
        Scheme & Memory consumption (\si{\mega \byte}) \\
        \hline
        Explicit & \num{25} \\
        Newton-Raphson & \num{127} \\
        $\nbLevelMax = 1$ & \num{63} \\
        $\nbLevelMax = 2$ & \num{74} \\
        $\nbLevelMax = 3$ & \num{77} \\
        $\nbLevelMax = 4$ & \num{78} \\
    \end{tabular}
    \caption[Memory consumption for the beam problem using the Jacobi method and the GMG algorithm]{Memory consumption for the explicit solver, the implicit solver using Newton-Raphson method and the geometric multigrid algorithm different values for $\nbLevelMax$.}
    \label{tab:gmg:beam:memory}
\end{table}

\Cref{tab:gmg:beam:memory} shows the memory consumption with the different methods.
The Newton-Raphson method requires solving large sparse linear systems.
For performance reasons, we store the matrix and its preconditioner.
This leads to a higher memory footprint than the other methods, \SI{127}{\mega \byte} instead of \SI{25}{\mega \byte} for the explicit solver.
Furthermore, the choice of the preconditioner has an impact on the memory consumption of the method.
This is discussed in \cite{bloch2021}.
We use a ``matrix-free'' approach for the Jacobi method, i.e., we do not store the operator $\jacobiOperator$, but we access it by computing $\jacobiOperator \left( \vector{\jacobiUnknown} \right)$. 
However, using the geometric multigrid algorithm requires storing temporary values at coarse levels, hence the increase of memory consumption with $\nbLevelMax$, from \SI{63}{\mega \byte} with $\nbLevelMax = 1$ to \SI{78}{\mega \byte} with $\nbLevelMax = 4$.
As $\nbLevelMax$ increases, more values have to be stored, but each level is coarser than the previous one, therefore fewer variables per additional level are needed.
For example, with $257 \times 257$ cells at the fine level and $\nbLevelMax = 4$, there are only $32 \times 32$ additional cells compared to $\nbLevelMax = 3$.
The code used here was developed to show the interest of the method, it was not optimized to minimize memory footprint.
For the sake of simplicity, we store at each level $\vector{\jacobiRhs}$, $\vector{\jacobiUnknown}$, $\vector{\jacobiIntermediate}$, $\widetilde{\vector{\jacobiUnknown}}$, $\vector{\jacobiError}$ and some auxiliary vectors at the coarsest levels $\coarseOperator \left( \coarseVector{\jacobiUnknown} \right)$, $\restrictionOperator{\spaceStep}{2\spaceStep}{\fineOperator \left( \fineVector{\jacobiUnknown} \right)}$, and $\coarseOperator \left( \coarseVector{\jacobiUnknown} \right) -\restrictionOperator{\spaceStep}{2\spaceStep}{\fineOperator \left( \fineVector{\jacobiUnknown} \right)}$.

Performance results obtained here are obtained with small configurations: $257 \times 257$ cells for a two-dimensional problem.
One can expect a better speed-up in the three-dimensional case.

\subsection{2D Riemann problem} \label{sect:gmg:riemannPb}

\begin{figure}
    \centering
    \newcommand{\multiplicationFactor}{5.3}

\begin{tikzpicture}
  \draw[dashed] (\multiplicationFactor*0, \multiplicationFactor*0.5) -- (\multiplicationFactor*1, \multiplicationFactor*0.5);
  \draw[dashed] (\multiplicationFactor*0.5, \multiplicationFactor*0) -- (\multiplicationFactor*0.5, \multiplicationFactor*1);

  \draw (\multiplicationFactor*0, \multiplicationFactor*0) -- (\multiplicationFactor*0, \multiplicationFactor*1);
  \draw (\multiplicationFactor*0, \multiplicationFactor*0) -- (\multiplicationFactor*1, \multiplicationFactor*0);
  \draw (\multiplicationFactor*0, \multiplicationFactor*1) -- (\multiplicationFactor*1, \multiplicationFactor*1);
  \draw (\multiplicationFactor*1, \multiplicationFactor*0) -- (\multiplicationFactor*1, \multiplicationFactor*1);

  \node[below] at (\multiplicationFactor*0, \multiplicationFactor*0) {\num{0}};
  \node[below] at (\multiplicationFactor*0.5, \multiplicationFactor*0) {\num{0.5}};
  \node[below] at (\multiplicationFactor*1, \multiplicationFactor*0) {\num{1}};

  \node[left] at (\multiplicationFactor*0, \multiplicationFactor*0) {\num{0}};
  \node[left] at (\multiplicationFactor*0, \multiplicationFactor*0.5) {\num{0.5}};
  \node[left] at (\multiplicationFactor*0, \multiplicationFactor*1) {\num{1}};

  \node at (\multiplicationFactor*0.25, \multiplicationFactor*0.25) {$\vector{\hat{F}} = 
    \begin{pmatrix}
      1\\0
    \end{pmatrix}
  \rightarrow$};

  \node at (\multiplicationFactor*0.75, \multiplicationFactor*0.25) {$\vector{\hat{F}} = 
    \begin{pmatrix}
      0\\-1
    \end{pmatrix}
  \downarrow$};

  \node at (\multiplicationFactor*0.25, \multiplicationFactor*0.75) {$\vector{\hat{F}} = 
    \begin{pmatrix}
      0\\1
    \end{pmatrix}
  \uparrow$};

  \node at (\multiplicationFactor*0.75, \multiplicationFactor*0.75) {$\vector{\hat{F}} = 
    \begin{pmatrix}
      -1\\0
    \end{pmatrix}
  \leftarrow$};

\end{tikzpicture}
    \caption[Initial condition for 2D Riemann problem]{Initial condition for 2D Riemann problem.}
    \label{fig:gmg:riemannPb_ci}
\end{figure}
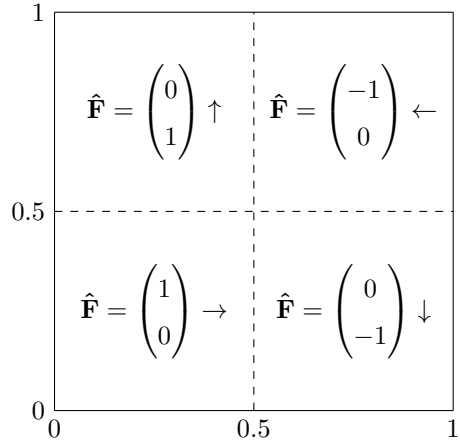

Let us now consider the same two-dimensional Riemann problem as in \cite{blachere2016}.
The domain $[0,1] \times [0,1]$ is discretized with $257 \times 257$ cells.
The initial temperature is $\gasTemperatureScheme_0 = \radiativeTemperature = \SI{1000}{\kelvin}$.
The domain is cut into four states, in each of them, the initial radiative flux is constant.
It is set to $\left( 1 - 10^{-8} \right) \speedOfLight \radiativeEnergy \vector{\hat{F}}$, with $\vector{\hat{F}}$ given by \cref{fig:gmg:riemannPb_ci}.

\begin{figure}
    \centering
    \begin{tikzpicture}
  \begin{semilogyaxis}[
      xlabel = {Number of time steps},
      ylabel = {Relative error on radiative energy},
      axis lines = left,
    ]
    \addplot[color=blue] table {images/riemannPb/output_energy_hll.txt};
    \addlegendentry{Explicit};
    
    \addplot[color=red] table {images/riemannPb/output_energy_jacobi.txt};
    \addlegendentry{$\nbLevelMax = \num{1}$};

    \addplot[color=green] table {images/riemannPb/output_energy_level2.txt};
    \addlegendentry{$\nbLevelMax = \num{2}$};
    
  \end{semilogyaxis}
\end{tikzpicture}
    \caption[Relative error for the radiative energy for the two-dimensional Riemann problem]{
    Evolution of the relative error on radiative energy as a function of time.}
    \label{fig:gmg:riemannPb_energy}
\end{figure}
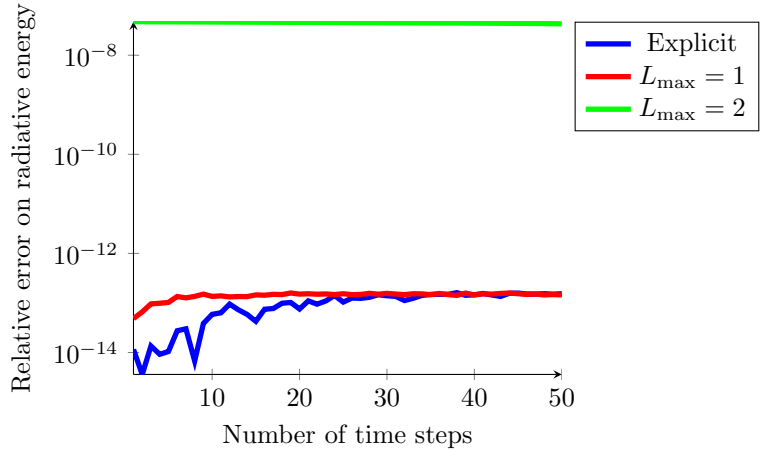

Using periodic boundary conditions, no energy should enter nor leave the box, therefore, the total radiative energy should be conserved at the precision
of the solver.
\Cref{fig:gmg:riemannPb_energy} shows the evolution of the relative error between the expected total radiative energy and the one actually computed in the box at each time step.
Using both the explicit solver and the Jacobi method, with the same time step $\Delta t$, the relative error
stays below \num{e-12}.
Using the GMG technique, the relative error is around \num{e-8}. 
The residual $\frac{\norm{\cycleIter{\fineVector{\jacobiResidual}}{\indexNbCycle}}}{\norm{\cycleIter{\fineVector{\jacobiResidual}}{0}}}$ is set to \num{5e-3}. 
Even for a quite high value of residual, the scheme is conservative.

\begin{figure}
    \centering
    \begin{tikzpicture}
  \begin{semilogyaxis}[
      xlabel = {Number of V-cycles ($\indexNbCycle$)},
      ylabel = {Residual $\frac{\norm{\cycleIter{\fineVector{\jacobiResidual}}{\indexNbCycle}}}{\norm{\cycleIter{\fineVector{\jacobiResidual}}{0}}}$},
      axis lines = left,
    ]
    \addplot[color=blue] table {images/riemannPb/output_level1.txt};
    \addlegendentry{$\nbLevelMax = 1$};

    \addplot[color=red] table {images/riemannPb/output_level2.txt};
    \addlegendentry{$\nbLevelMax = 2$};

    \addplot[color=green] table {images/riemannPb/output_level3.txt};
    \addlegendentry{$\nbLevelMax = 3$};

    \addplot[color=Orange] table {images/riemannPb/output_level4.txt};
    \addlegendentry{$\nbLevelMax = 4$};

  \end{semilogyaxis}
\end{tikzpicture}
    \caption[Residual as a function of the number of V-cycle for the two-dimensional Riemann problem]{Evolution of the residual as a function of the number of V-cycle, with different values for $\nbLevelMax$ for the two-dimensional Riemann problem.}
    \label{fig:gmg:riemannPb_cv}
\end{figure}
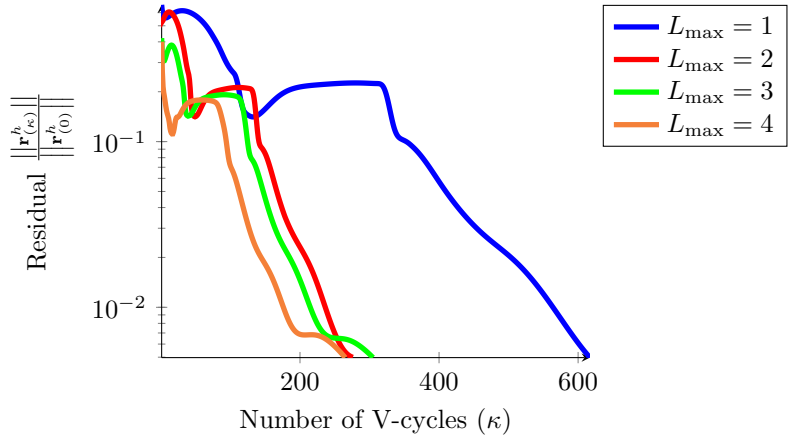

As shown by \cref{fig:gmg:riemannPb_cv}, with a large time step ($\cfl = \num{2000}$), increasing the value of $\nbLevelMax$ does not help to reduce the number of V-cycles to reach the same residual.
When increasing the value of $\nbLevelMax$, the computational time per V-cycle increases, therefore the total computational time also increases, unlike the beam test.
Because we are in the free-streaming regime, the propagation of the photons has to be followed.
It can be done either by reducing the time step or by performing numerous iterations of the Jacobi method with a large time step.
We use the latter to highlight the properties of the method.

As in \cref{sect:gmg:beam}, all the curves have the same shape with different values of $\nbLevelMax$, but it differs from the curves obtained with the beam test.

\begin{figure}
    \begin{subfigure}[b]{0.5\textwidth}
        \begin{centering}
            \includegraphics[width=\columnwidth]{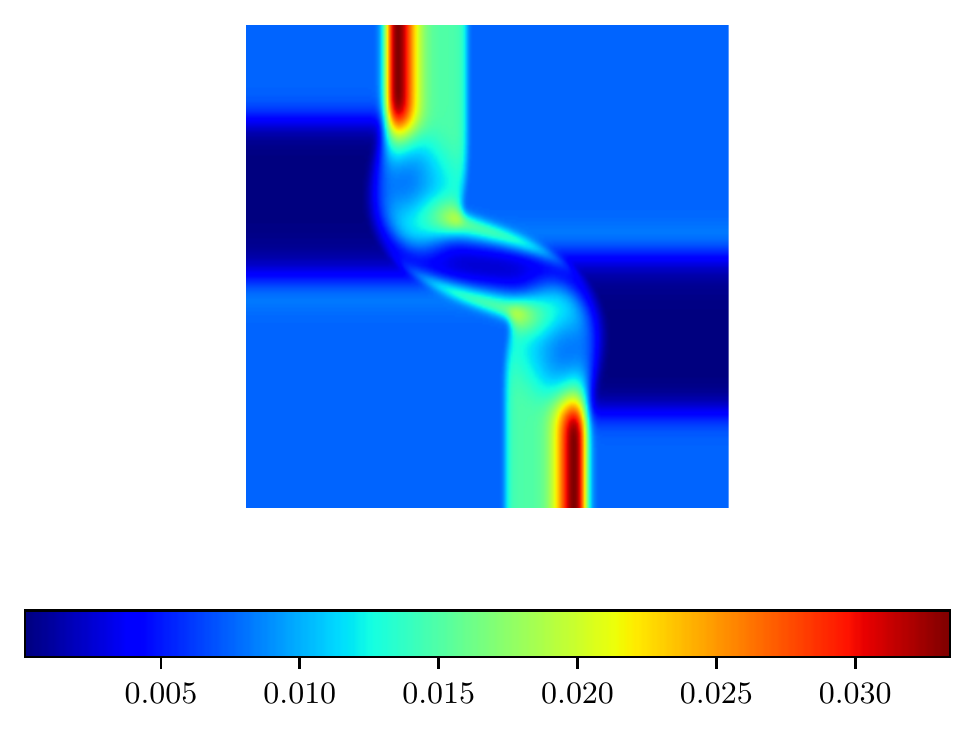}
            \caption{Explicit, $\cfl = \num{0.9}$}
            \label{fig:gmg:riemannPb:hll}
        \end{centering}
    \end{subfigure}
    \hfill
    \begin{subfigure}[b]{0.5\textwidth}
        \begin{centering}
            \includegraphics[width=\columnwidth]{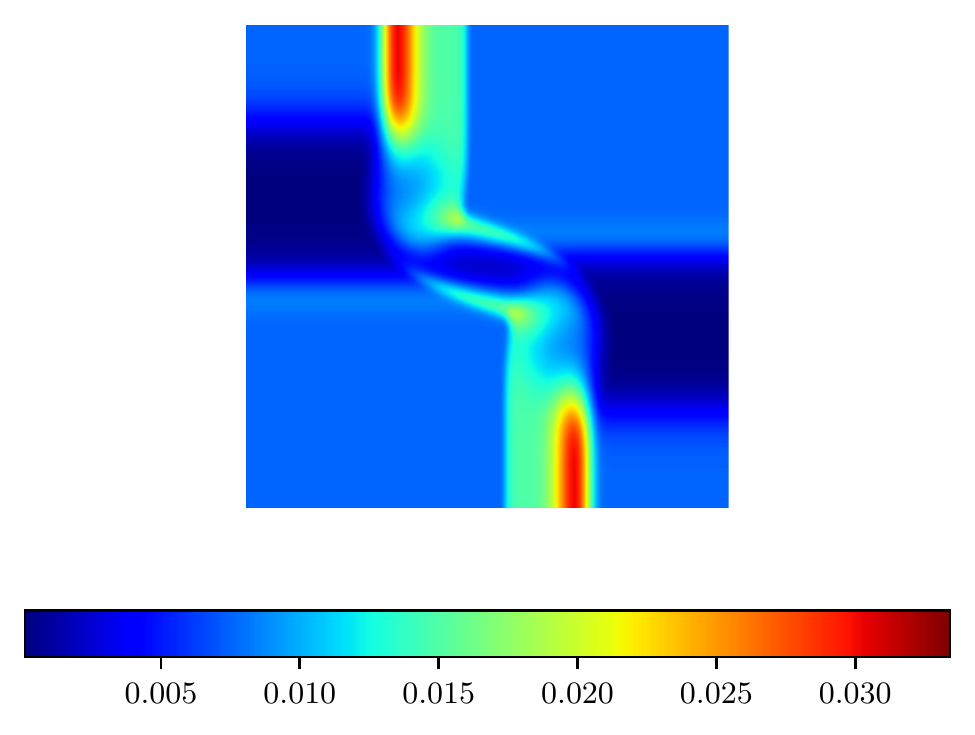}
            \caption{Jacobi method, $\cfl = \num{0.9}$}
            \label{fig:gmg:riemannPb:explicitTimeStep}
        \end{centering}
    \end{subfigure}
%
    \begin{subfigure}[b]{0.5\textwidth}
        \begin{centering}
            \includegraphics[width=\columnwidth]{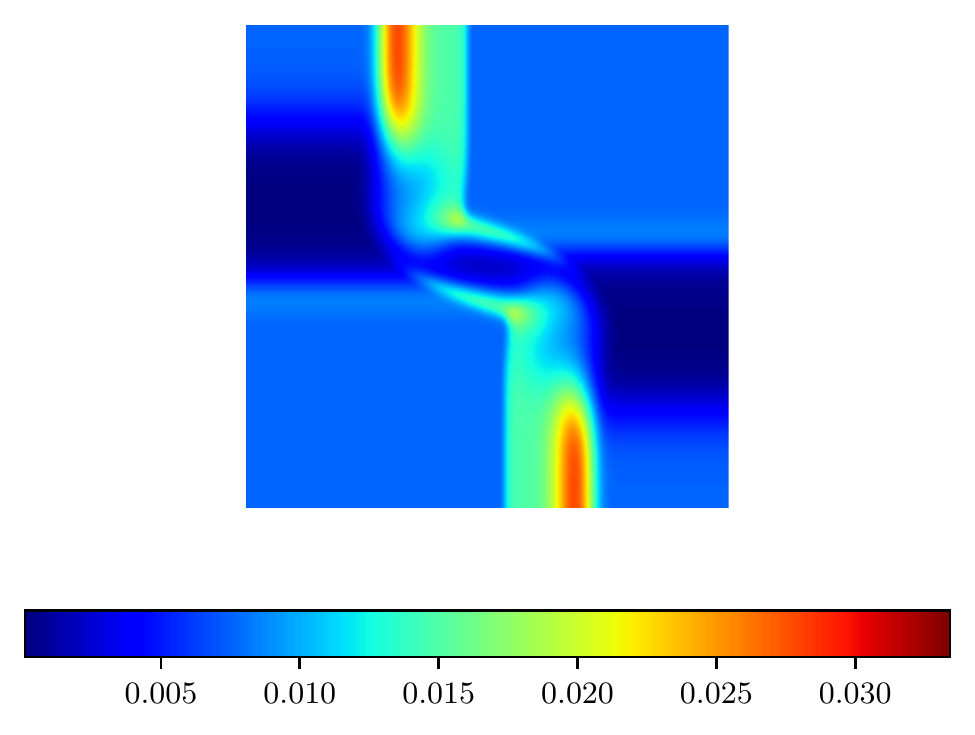}
            \caption{Jacobi method, $\cfl = \num{2.7}$}
            \label{fig:gmg:riemannPb:smallTimeStep}
        \end{centering}
    \end{subfigure}
    \hfill
    \begin{subfigure}[b]{0.5\textwidth}
        \begin{centering}
            \includegraphics[width=\columnwidth]{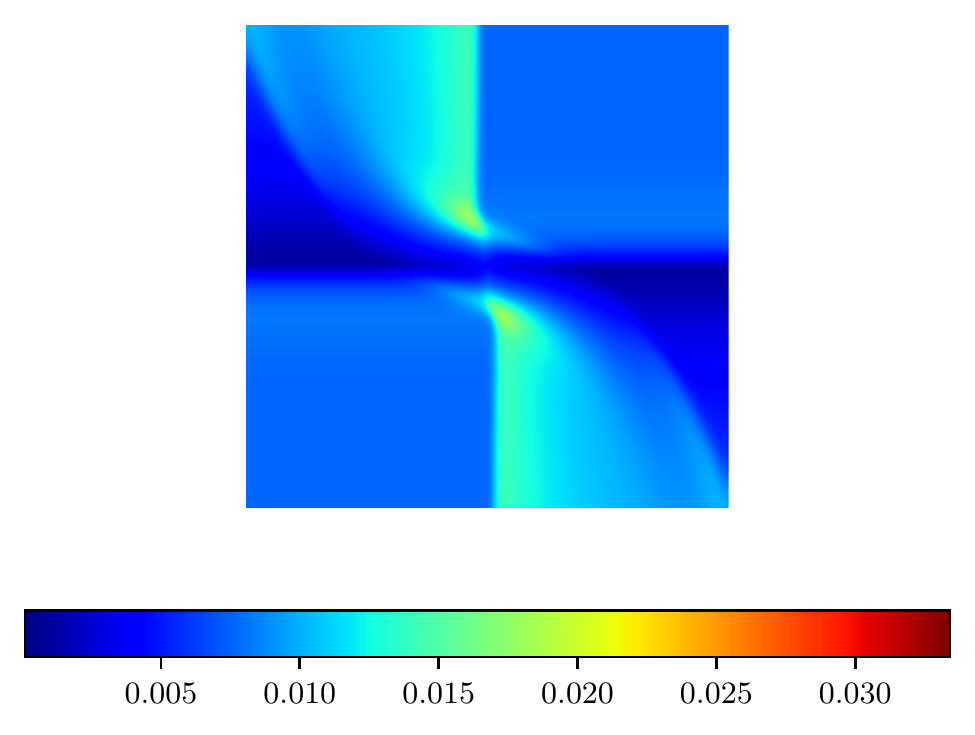}
            \caption{Jacobi method, $\cfl = \num{150}$}
            \label{fig:gmg:riemannPb:largeTimeStep}
        \end{centering}
    \end{subfigure}
    \caption[Radiative energy in two-dimensional Riemann problem with the Jacobi method]{Snapshots of radiative energy at final time $t_f = \SI{e-11}{\second}$ with the explicit solver and the Jacobi method with different time steps.}
    \label{fig:gmg:riemannPb}
\end{figure}

\Cref{fig:gmg:riemannPb} shows snapshots of radiative energy at the final time \SI{e-11}{\second} when Neumann boundary conditions are used.
Different solutions are obtained with the explicit solver (\cref{fig:gmg:riemannPb:hll}) and Jacobi method with different time steps (\cref{fig:gmg:riemannPb:explicitTimeStep,fig:gmg:riemannPb:smallTimeStep,fig:gmg:riemannPb:largeTimeStep}).
The residual for the Jacobi method is set to \num{1e-6}.
As discussed in \cite{blachere2016}, the reduced flux stays close to \num{1} during the simulation and the method is stable even for very stiff problems, with any time step.
Because we use a scheme of order \num{1}, the solution is diffused, but the method preserves the admissible states with the reduced flux close to \num{1}.
When small time steps are used with the Jacobi method (\cref{fig:gmg:riemannPb:explicitTimeStep,fig:gmg:riemannPb:smallTimeStep}), we recover the explicit solution.
When larger time steps are used (\cref{fig:gmg:riemannPb:largeTimeStep}), the solution is more diffused. 
Because the flow is in the free-streaming regime and with sharp discontinuities, the propagation of the photons should be followed to capture with precision the solution of this problem.
This can only be done by reducing the time step, close to the explicit time step.
We do not expect an implicit method to be efficient in the test, but it shows that our scheme is extremely stable and robust, even at very large time steps.
This feature is essential to carry out large and complex radiation hydrodynamics simulations, where the physics of interest will mostly be on hydrodynamical time scales.

\section{Discussion and conclusion} \label{sect:gmg:conc}

We have first presented a Jacobi method to solve the nonlinear system arising from the discretization of the M\textsubscript{1} model without source terms with a time-implicit HLL solver. 
The admissible states are preserved, even with large CFL numbers.
This method is iterative, the convergence rate decreases when the resolution increases.
To tackle this issue, we use a nonlinear geometric multigrid algorithm.
The Jacobi method first described is used as smoother and coarse grid solver.
However, this algorithm relies on the residual equation, which does not preserve the admissible states.
Instead of solving this equation, we introduce a pseudo-time, and we look for a steady state in pseudo-time.
Numerical experiments have shown the good performances obtained with this algorithm.
However, it should be used carefully in physical problems that require following the propagation of photons.

\subsection{Source terms}

In order to study physical problems, the integration of source terms is indispensable.
It needs to be done carefully, regarding the asymptotic preserving property discussed in \cite{gonzalez2007,bloch2021} for example.
The scheme should capture an asymptotic behavior in the diffusive limit.
We did not address this issue here.
For example,
schemes presented in \cite{jin1996,buetDespres2008,berthon2011} are explicit, they preserve the admissible states, and they are
asymptotic preserving.
The HLL solver used here could be replaced by
another
scheme.
This will change the nonlinear operator $\jacobiOperator$.
It requires choosing carefully which terms are taken at iteration $\jacobiCounter$ and $\jacobiCounter+1$ in \cref{algo:jacobi}.
The geometric multigrid algorithm presented in \cref{sect:gmg:multigrille} does not rely on the form of the nonlinear operator, as long as the iterations of the Jacobi method preserve the admissible states.
Therefore, the geometric multigrid method can be used to study physical problems with opacity as soon as the Jacobi algorithm is extended to an asymptotic preserving scheme.

Similarly, the method could be extended to a non-gray model.
We have defined the radiative energy (resp. radiative flux and radiative pressure) as the mean over the frequency and the direction of propagation of the photons of the specific intensity $I$ (resp. $I \vector{\Omega}$ and $I \vector{\Omega} \otimes \vector{\Omega}$).
Therefore, our model does not take into account the frequency of the photons.
One can use a multigroup model.
The continuous frequency is discretized into several groups, and the variables are integrated over the frequency in each group.
As previously, the Jacobi method could be adapted to this model. 
The main difficulty is to choose which terms are taken at iteration $\jacobiCounter$ and $\jacobiCounter+1$ in \cref{algo:jacobi}.

\subsection{High-order scheme}

In the numerical examples presented in \cref{sect:gmg:resultatsNum}, the solution is more diffused than expected because a scheme of order \num{1} is used.
The numerical solution can be improved by increasing the order of the scheme.
One method, among others, is spectral volume method (see e.g., \cite{padioleau2020}).
Several degrees of freedom are added in the cell to reconstruct a polynomial.
As with the finite volume method, the solution is discontinuous at the interfaces of the cells to capture shocks.
This method is different from the Discontinuous Galerking method because it is based on the integral formulation of the equation instead of the weak formulation.
The Runge-Kutta method can be used to achieve high-order in time.

Because the spectral volume method is close to the finite volume method, the Jacobi algorithm can be adapted to this scheme.
As mentioned in \cite{kifonidis2012}, the pseudo-time introduced in \cref{sect:gmg:m1} is independent of the physical time.
The discretization of \cref{eq:gmg:pseudoTime:pseudoTimeEqua} is not required to be accurate, and an integration of order \num{1} could be used.

In addition, a $p$-multigrid method can be used.
A hierarchy is build and each level corresponds to a different approximation order.

However, we think that a high-order version of the algorithm is a significant amount of work, especially the inclusion of the reconstrution step in the non-linear jacobi and GMG method.
This is part of our perspective for the future.

\subsection{Performances}

Although the performance results shown in \cref{sect:gmg:resultatsNum} are promising, some choices were made and others were not explored.
For example, only a V-cycle is used, but there exist other possibilities: W-cycle, F-cycle,\ldots
The restriction and prolongation operators can also have an impact on the performances.

Another well-known method to reduce the number of iterations performed by the Jacobi method is to replace it with a Gauss-Seidel algorithm.
Even though the tests presented in \cref{sect:gmg:resultatsNum} are obtained with a sequential code, we aim at using a parallel implementation to study physical situations.
In the linear case, the Jacobi method is known to be easier to parallelize than the Gauss-Seidel method. 

\subsection{Parallel computing}

The operator $\jacobiOperator$ considered here is nonlinear, due to the nonlinearity of the model.
However, if this operator were linear, \cref{algo:jacobi} would be the classical Jacobi method for a matrix-free problem.
The pros and cons are the same as any multigrid method.
For realistic computations, parallel computing is mandatory.
The Jacobi method itself is quite easy to parallelize, as in the linear case.
Because our approach can be seen as a matrix-free method, it can be parallelized using techniques developed for other equations.
See for example \cite{rude2021} about the resolution of the Poisson equation using a parallel matrix-free multigrid method.
Matrix-free algorithms are very suited for massively parallel computing because of their low memory footprint.

\section*{Acknowledgements}

PT acknowledges supports by the European Research Council under Grant Agreement ATMO 757858.

\appendix
\section{Overall algorithm} \label{appendix:gmg:overallAlgo}

\Cref{algo:gmgOverall} presents the whole algorithm to compute $\vector{\jacobiUnknown}^{n+1}$ with $\vector{\jacobiUnknown}^{n}$ known.

\begin{breakablealgorithm}
    \caption{Overall algorithm}
    \begin{algorithmic}[1]
        \State Choose a maximum number of levels $\nbLevelMax$
        \State $\indexNbCycle = 0$
        \State $\fineVector{\jacobiRhs} = \vector{\jacobiUnknown}^n$
        \State $\cycleIter{\fineVector{\jacobiUnknown}}{0} = \vector{\jacobiUnknown}^n$. From now on, the exponent is the mesh size, it is no longer the time
        \While{$\norm{\cycleIter{\vector{\jacobiResidual}}{\indexNbCycle}} > \tolGmg$}
        \For{$\indexLevel = 0, \cdots, \nbLevelMax-1$}
        \If{$\indexLevel = 0$}
        \State Pre-smoother: solve $\fineOperator \left( \cycleIter{\fineVector{\jacobiIntermediate}}{\indexNbCycle} \right) = \fineVector{\jacobiRhs}$ with $\smootherNbIter{0}$ iterations of the nonlinear Jacobi method (\cref{algo:jacobi}), with initial guess $\cycleIter{\fineVector{\jacobiUnknown}}{\indexNbCycle}$
        \State Compute residual: $\cycleIter{\fineVector{\overline{\jacobiResidual}}}{\indexNbCycle} = \fineVector{\jacobiRhs} - \fineOperator \left( \cycleIter{\fineVector{\jacobiIntermediate}}{\indexNbCycle} \right)$
    \Else
        \State Pre-smoother: compute the steady state in pseudo-time of $\frac{\dif \cycleIter{\vectorAtLevel{\jacobiIntermediate}{\indexLevel}}{\indexNbCycle}}{\dif \pseudoTime} + \operatorAtLevel{\indexLevel} \left( \cycleIter{\vectorAtLevel{\jacobiIntermediate}{\indexLevel}}{\indexNbCycle} \right) = \operatorAtLevel{\indexLevel} \left( \cycleIter{\vectorAtLevel{\jacobiUnknown}{\indexLevel}}{\indexNbCycle} \right) + \cycleIter{\vectorAtLevel{\jacobiResidual}{\indexLevel}}{\indexNbCycle}$ using \cref{algo:pseudoTimeTwoGrids}, with $\smootherNbIter{\indexLevel}$ iterations of the nonlinear Jacobi method , with initial guess $\cycleIter{\vectorAtLevel{\jacobiUnknown}{\indexLevel}}{\indexNbCycle}$
        \State Compute residual: $\cycleIter{\vectorAtLevel{\overline{\jacobiResidual}}{\indexLevel}}{\indexNbCycle} = \operatorAtLevel{\indexLevel} \left( \cycleIter{\vectorAtLevel{\jacobiUnknown}{\indexLevel}}{\indexNbCycle} \right) + \cycleIter{\vectorAtLevel{\jacobiResidual}{\indexLevel}}{\indexNbCycle} - \operatorAtLevel{\indexLevel} \left( \cycleIter{\vectorAtLevel{\jacobiIntermediate}{\indexLevel}}{\indexNbCycle} \right)$
        \EndIf
        \State Restriction: $\cycleIter{\vectorAtLevel{\jacobiUnknown}{\indexLevel+1}}{\indexNbCycle} = \restrictionOperator{2^{\indexLevel} \spaceStep}{2^{\indexLevel+1} \spaceStep}{\cycleIter{\vectorAtLevel{\jacobiUnknown}{\indexLevel}}{\indexNbCycle}}$ and $\cycleIter{\vectorAtLevel{\jacobiResidual}{\indexLevel+1}}{\indexNbCycle} = \restrictionOperator{2^{\indexLevel} \spaceStep}{2^{\indexLevel+1} \spaceStep}{\cycleIter{\vectorAtLevel{\overline{\jacobiResidual}}{\indexLevel}}{\indexNbCycle}}$
        \EndFor
        \State Solver at the coarsest level: compute the steady state in pseudo-time of $\frac{\dif \cycleIter{\vectorAtLevel{\jacobiIntermediate}{\nbLevelMax}}{\indexNbCycle}}{\dif \pseudoTime} + \operatorAtLevel{\nbLevelMax} \left( \cycleIter{\vectorAtLevel{\jacobiIntermediate}{\nbLevelMax}}{\indexNbCycle} \right) = \operatorAtLevel{\nbLevelMax} \left( \cycleIter{\vectorAtLevel{\jacobiUnknown}{\nbLevelMax}}{\indexNbCycle} \right) + \cycleIter{\vectorAtLevel{\jacobiResidual}{\nbLevelMax}}{\indexNbCycle}$ using \cref{algo:pseudoTimeTwoGrids}, with $\smootherNbIter{\nbLevelMax}$ iterations of the nonlinear Jacobi method , with initial guess $\cycleIter{\vectorAtLevel{\jacobiUnknown}{\nbLevelMax}}{\indexNbCycle}$
        \For{$\indexLevel = \nbLevelMax-2, \cdots, 0$ by step $-1$}
        \State Prolongation: $\cycleIter{\vectorAtLevel{\overline{\jacobiUnknown}}{\indexLevel}}{\indexNbCycle} = \cycleIter{\vectorAtLevel{\jacobiUnknown}{\indexLevel}}{\indexNbCycle} + \prolongationOperator{2^\indexLevel \spaceStep}{2^{\indexLevel+1} \spaceStep}{\cycleIter{\vectorAtLevel{\jacobiIntermediate}{\indexLevel+1}}{\indexNbCycle} - \cycleIter{\vectorAtLevel{\jacobiUnknown}{\indexLevel+1}}{\indexNbCycle}}$ 
        \For{each cell $\indexCellLevel{i}{\indexLevel}$ in the grid at level $\indexLevel$}
        \If{$\left( \cycleIter{\vectorAtLevel{\overline{\jacobiUnknown}}{\indexLevel}}{\indexNbCycle} \right)_{\indexCellLevel{i}{\indexLevel}}$ is not admissible}
        \State use \cref{eq:gmg:pseudoTime:prolongatedSol} to update $\left( \cycleIter{\vectorAtLevel{\overline{\jacobiUnknown}}{\indexLevel}}{\indexNbCycle} \right)_{\indexCellLevel{i}{\indexLevel}}$
        \EndIf
        \EndFor
        \If{$\indexLevel = 0$}
        \State Post-smoother: solve $\fineOperator \left( \cycleIter{\fineVector{\jacobiIntermediate}}{\indexNbCycle} \right) = \fineVector{\jacobiRhs}$ with $\smootherNbIter{0}$ iterations of the nonlinear Jacobi method (\cref{algo:jacobi}), with initial guess $\cycleIter{\fineVector{\overline{\jacobiUnknown}}}{\indexNbCycle}$
    \Else
        \State Post-smoother: compute the steady state in pseudo-time of $\frac{\dif \cycleIter{\vectorAtLevel{\jacobiIntermediate}{\indexLevel}}{\indexNbCycle}}{\dif \pseudoTime} + \operatorAtLevel{\indexLevel} \left( \cycleIter{\vectorAtLevel{\jacobiIntermediate}{\indexLevel}}{\indexNbCycle} \right) = \operatorAtLevel{\indexLevel} \left( \cycleIter{\vectorAtLevel{\jacobiUnknown}{\indexLevel}}{\indexNbCycle} \right) + \cycleIter{\vectorAtLevel{\jacobiResidual}{\indexLevel}}{\indexNbCycle}$ using \cref{algo:pseudoTimeTwoGrids}, with $\smootherNbIter{\indexLevel}$ iterations of the nonlinear Jacobi method , with initial guess $\cycleIter{\vectorAtLevel{\overline{\jacobiUnknown}}{\indexLevel}}{\indexNbCycle}$
        \EndIf
        \EndFor
        \State $\cycleIter{\fineVector{\jacobiUnknown}}{\indexNbCycle} = \cycleIter{\fineVector{\jacobiIntermediate}}{\indexNbCycle}$
        \State $\cycleIter{\vector{\jacobiResidual}}{\indexNbCycle} = \vector{\jacobiUnknown}^n - \fineOperator \left( \cycleIter{\fineVector{\jacobiUnknown}}{\indexNbCycle} \right)$
        \If{$\frac{\abs{ \norm{\cycleIter{\vector{\jacobiResidual}}{\indexNbCycle}} - \norm{\cycleIter{\vector{\jacobiResidual}}{\indexNbCycle-1}} }}{\norm{\cycleIter{\vector{\jacobiResidual}}{0}}} > \tolPseudoTimeIncrease$}
        \State $\dTauIm \gets 1.1 \dTauIm$
        \EndIf
        \If{$\frac{\abs{ \norm{\cycleIter{\vector{\jacobiResidual}}{\indexNbCycle}} - \norm{\cycleIter{\vector{\jacobiResidual}}{\indexNbCycle-1}} }}{\norm{\cycleIter{\vector{\jacobiResidual}}{0}}} < \tolPseudoTimeDecrease$}
        \State $\dTauIm \gets \dTauIm / 2$
        \EndIf
        \State $\indexNbCycle \gets \indexNbCycle+1$
        \EndWhile
        \State $\vector{\jacobiUnknown}^{n+1} = \cycleIter{\fineVector{\jacobiUnknown}}{\indexNbCycle}$
    \end{algorithmic}
    \label{algo:gmgOverall}
\end{breakablealgorithm}

\bibliography{mybibfile}

\end{document}